\documentclass{article}

\usepackage[margin=1in]{geometry}
\usepackage{amsmath,amssymb,amsthm,mathtools}
\usepackage{microtype}
\usepackage{enumitem}
\usepackage[hidelinks]{hyperref}
\usepackage{xcolor}
\usepackage{tikz}

\usepackage[commentmarkup=uwave]{changes}
\definechangesauthor[name={CG},color=magenta]{cg}
\definechangesauthor[name={MP},color=red]{mp}
\definechangesauthor[name={MK},color=blue]{mk}
\definechangesauthor[name={AI},color=magenta]{ai}

\usepackage{lineno}

\usepackage{cleveref}

\newtheorem{theorem}{Theorem}[section]
\newtheorem{example}{Example}[section]
\newtheorem{lemma}[theorem]{Lemma}
\newtheorem{conjecture}{Conjecture}[section]
\newtheorem{proposition}[theorem]{Proposition}
\newtheorem{corollary}[theorem]{Corollary}

\newcommand{\PP}{\mathbb P}
\newcommand{\EE}{\mathbb E} 
\newcommand{\Indic}{\mathbf 1}
\newcommand{\N}{\mathbb N}

\title{Graphon hitting times: connectivity and minimum degree}
\author{Constanza Gacit\'ua Fuentes\thanks{Departamento de Ingenier\'ia Matemática and Centro de Modelamiento Matemático, Universidad de Chile, Santiago, Chile. \texttt{E-mail:\{cgacitua,mk,mpavez\}@dim.uchile.cl}. All the authors acknowledge the support of BASAL funds for centers of excellence from ANID-Chile (FB210005). MPS is also supported by ANID Fondecyt Regular grant No. 1241398.} \and Marcos Kiwi\footnotemark[1] \and Mat\'ias Pavez-Sign\'e\footnotemark[1]}
\date{}

\begin{document}
\maketitle

\begin{abstract}
A classical result in random graphs states that the Erd\H{o}s--R\'enyi random graph process becomes connected at exactly the moment its last isolated vertex disappears. The main goal of this paper is to establish a graphon analogue of this result. Let $G(n,W)$ be the inhomogeneous random graph sampled from a graphon $W$, and reveal its edges one by one in a uniformly random order. We prove that if $W$ is connected and the lower tail of its degree distribution satisfies $g_W(t)/t\to 0$ as $t\to 0$, which is needed for the connectivity of $G(n,W)$, then the connectivity and minimum-degree-one hitting times coincide asymptotically almost surely. We also determine when this hitting time occurs at scale $\Theta(n\log n)$, as in the binomial random graph model. Moreover, we give a simple analytic criterion for the transition to be localized in a $(1+o(1))$-multiplicative window, and construct graphons with minimum degrees arbitrarily close to zero whose transition is nevertheless sharp. Our approach uses Laplace transforms and vague convergence of rescaled graphon degree measures to provide a unified framework for analyzing low-degree vertices and connectivity transitions.
Along the way, we revisit the minimum degree of~$G(n,W)$, giving simplified proofs of results of Hladk\'y and Viswanathan~\cite{HV} and establishing the limiting distribution of $\delta(G(n,W))$ under a second-order decorrelation condition.
\end{abstract}

\section{Introduction}\label[section]{sec:intro}
The binomial random graph model $G(n,p)$ is the random graph with vertex set $[n]$ where each of the possible~$\binom{n}{2}$ edges appears independently with probability $p$. Since the seminal work by Erd\H os and R\'enyi~\cite{ErdosRenyi1959}, understanding the connectivity of $G(n,p)$ has become a central problem in the field. In fact, the result by Erd\H os and R\'enyi implies that, for any fixed~$\varepsilon>0$, the graph $G(n,p)$ is asymptotically almost surely\footnote{We say that a sequence of events $(A_n)_{n\in\mathbb N}$ holds asymptotically almost surely if $\lim_{n\to\infty} \mathbb P(A_n)=1$.} (a.a.s.~for short) connected when $p\ge (1+\varepsilon)\frac{\log n}{n}$, and a.a.s. disconnected if $p\le (1-\varepsilon)\frac{\log n}{n}$. Moreover, the thresholds for connectivity and the disappearance of isolated vertices in $G(n,p)$ coincide: for every fixed~$\varepsilon>0$, the graph $G(n,p)$ a.a.s.\ contains no isolated vertices when $p\ge (1+\varepsilon)\frac{\log n}{n}$, whereas it a.a.s.\ contains an isolated vertex when $p\le (1-\varepsilon)\frac{\log n}{n}$.

In this paper, we study the connectivity properties of inhomogeneous random graphs sampled from \emph{graphons}, where a {graphon} is a Lebesgue measurable symmetric function $W\colon [0,1]^2\to [0,1]$. Graphons are important objects in extremal combinatorics, as they arise as the limit object of dense graph sequences (see e.g.~\cite{LovaszSzegedy06,BCLSV08,BCLSV12,LovaszBook,BCL10,Janson13}), and, moreover, provide a broad class of inhomogeneous random graph models, as we define next. Given  $n\in\mathbb N$ and a graphon $W:[0,1]^2\to [0,1]$, the $W$-random graph $G(n,W)$ is the random graph with vertex set $[n]$ which is defined in two stages, the vertex-generating stage and the edge-generating stage. First, in the vertex-generating stage, we sample~$n$ independent elements $X_1,\ldots, X_n\in [0,1]$ uniformly at random, and regard $X_i$ as the \emph{type} of the vertex~$i$. Second, in the edge-generating stage, each possible edge~$ij$ appears in $G$  with probability $W(X_i,X_j)$, independently of all other edges. Observe that the binomial random graph model $G(n,p)$ has the same distribution as $G(n,W)$ when $W\equiv p$ is the constant graphon.

In contrast to the binomial case, the random graph $G(n,W)$ may contain isolated vertices even when the resulting graph has $\Theta(n^2)$ edges. The main difference here is that $G(n,W)$ may be locally sparse even when~$G(n,W)$ is globally dense. For example, if $W(x,y)$ approaches 0 in some region $\Omega\subseteq[0,1]$, then $G(n,W)$ behaves similarly to a sparse random graph in those vertices with type in $\Omega$. In this case, the natural measure of local sparsity is the \emph{graphon degree function} $d_W:[0,1]\to [0,1]$, which is defined by
\[
d_W(x):=\int_0^1W(x,y)\,dy.\]
The degree function $d_W(\cdot)$ governs the degree of the vertices in $G(n,W)$, as a vertex of type $x\in [0,1]$ has (conditional) expected degree $(n-1)d_W(x)$. Indeed, Devroye and
Fraiman~\cite{DF} characterized the connectivity of random graphs sampled from {kernels} in terms of their degree function, which was recently revisited by Hladk\'y and Viswanathan~\cite{HV} in the case of graphons. To describe this connection more precisely, let~$S_W(t):=\{x\in[0,1]:d_W(x)\le t\}$ denote the set of low-degree types in $W$ and define the
\emph{lower-tail graphon degree distribution} of $W$ by
\[
g_W(t):=\lambda\big(S_W(t)\big),
\]
where $\lambda$ denotes the Lebesgue measure. Hladk\'y and Viswanathan showed that the behavior of the minimum
degree of $G(n,W)$ is controlled by the lower tail $g_W(t)$, as follows.
\begin{theorem}[Hladk\'y and Viswanathan~\cite{HV}]\label{thm:hladky} Let $W$ be a graphon.
\begin{enumerate}[label=(\alph{enumi})]
    \item If $\displaystyle\lim_{t\to 0}\frac{g_W(t)}{t}=+\infty$, then a.a.s.\ $G(n,W)$ contains at least one isolated vertex.
    \item  If $\displaystyle\lim_{t\to 0}\frac{g_W(t)}{t}=0$, then $\delta(G(n,W))\to +\infty$ in probability. 
    \item  If $\displaystyle\liminf_{t\to 0}\frac{g_W(t)}{t}>0$, then $\liminf_{n\to\infty}\mathbb P(G(n,W)\text{ contains an isolated vertex})>0$.
    \item If $\displaystyle\limsup_{t\to 0}\frac{g_W(t)}{t}<+\infty$, then for every $\ell\in\mathbb N$ we have $\liminf_{n\to\infty}\mathbb P(\delta(G(n,W))>\ell)>0$.
\end{enumerate}
\end{theorem}
Intuitively, a vertex of type $x\in [0,1]$ is likely to be isolated if $d_W(x)=o(1/n)$. 
The expected number of vertices with type-degree at most $t$ in $G(n,W)$ is $ng_W(t)$. Thus, for small $c>0$, the quantity $ng_W(c/n)$ controls the number of isolated vertices.

\subsection{The random graph process and hitting times}
The \emph{random graph process} $(G_t)_{t=0}^{N}$, with $N:=\binom{n}{2}$, starts from the empty graph $G_0$ on vertex set $[n]$, and at each step, $G_t$ is obtained from $G_{t-1}$ by adding a single edge chosen uniformly at random among the non-edges of $G_{t-1}$. Equivalently, the process is determined by a uniformly random ordering of the edges of $K_n$, and for each fixed $m$ the graph $G_m$ consists of the first $m$ edges in this ordering. This process formulation turns the threshold statements into much sharper \emph{hitting time} results. For a monotone increasing property $\mathcal{P}$, let $\tau(\mathcal{P}):=\min\{t : G_t\in\mathcal{P}\}$, and write $\tau_k:=\tau(\delta(G)\ge k)$ and $\kappa_k:=\tau(k\text{-connected})$. Since $k$-connectivity implies minimum degree at least $k$, we always have $\tau_k\le\kappa_k$. Bollob\'as and Thomason~\cite{BollobasThomason1985} proved that  $\tau_k=\kappa_k$ a.a.s.~for every fixed $k\ge1$, that is, the very edge that raises the minimum degree to $k$ also makes the graph~$k$-connected. In particular, for $k=1$, this says that the random graph becomes connected exactly when its last isolated vertex disappears. Results of the same flavor hold for other spanning properties whose only obstructions are local. Most notably, Bollob\'as~\cite{Bollobas1984} and, independently, Ajtai, Koml\'os and Szemer\'edi~\cite{AjtaiKomlosSzemeredi1985} showed that a.a.s.\ $G_t$ becomes Hamiltonian exactly at time $\tau_2$. 

A natural generalization of the random graph process is replacing the complete graph $K_n$ by an arbitrary fixed host graph $\Gamma$ on $n$ vertices. The \emph{random subgraph process} on $\Gamma$ is the sequence $(\Gamma_t)_{t=0}^{e(\Gamma)}$ obtained by fixing a uniformly random ordering of $E(\Gamma)$ and revealing its edges one by one, starting from the empty graph on $V(\Gamma)$ and finishing with $\Gamma$. Each $\Gamma_t$ is thus a uniformly random subgraph of $\Gamma$ with exactly $t$ edges. For connectivity, the minimum degree remains a trivial obstruction, so one again compares the hitting time~$\tau_k(\Gamma):=\min\{t:\delta(\Gamma_t)\ge k\}$ with $\tau_{k-\mathrm{conn}}(\Gamma):=\min\{t:\Gamma_t\text{ is $k$-connected}\}$. A closely related model is the random subgraph $\Gamma_p$, which is obtained by keeping each edge $e\in E(\Gamma)$ with probability $p$, making all choices independently.  In contrast with the case~$\Gamma=K_n$, the hitting times $\tau_k(\Gamma)$ and $\tau_{k-\mathrm{conn}}(\Gamma)$ do not need to coincide in general. For instance, if $\Gamma$ consists of two disjoint cliques joined by a single edge, then~$\Gamma_t$ is disconnected until that edge is revealed, which typically happens long after $\tau_1(\Gamma)$. In order to get that~$\tau_k(\Gamma)=\tau_{k-\mathrm{conn}}(\Gamma)$ a.a.s., one typically needs some sort of expansion or large minimum degree in the host graph.

Hitting-time phenomena for other spanning properties have also been extensively studied in random subgraph processes on structured host graphs.
One of the classical examples is the $n$-dimensional hypercube $Q_n$, for which Bollob\'as~\cite{Bollobas1990} proved that, in the random subgraph process on $\Gamma=Q_n$, a.a.s.\ a perfect matching appears exactly when the last isolated vertex disappears, and more recently, Condon, Espuny D\'iaz, Gir\~ao, K\"uhn and Osthus~\cite{CEGKO2024} showed that a.a.s.~the random subgraph process on $Q^n$ contains $k$ edge-disjoint Hamilton cycles as soon as its minimum degree reaches $2k$. Besides hypercubes, host graphs with large minimum degree have also received considerable attention. For instance,
Krivelevich, Lee and Sudakov~\cite{KrivelevichLeeSudakov2014} proved that there is a constant $C>0$
such that if $\Gamma$ is a Dirac graph, that is, $\delta(\Gamma)\geq n/2$, then the random subgraph~$\Gamma_p$ is a.a.s.\
Hamiltonian whenever $p\geq C\log n/n$. The corresponding hitting-time result was obtained by
Johansson~\cite{Johansson2020}, that is, if $\delta(\Gamma)\geq(1/2+\varepsilon)n$ for some fixed
$\varepsilon>0$, then a.a.s.\ the random subgraph process on $\Gamma$ becomes Hamiltonian exactly at
time $\tau_2(\Gamma)$. Alon and Krivelevich~\cite{AlonKrivelevich2022} extended this to~$k$
edge-disjoint Hamilton cycles, which a.a.s.\ appear exactly at time $\tau_{2k}(\Gamma)$, both under the
same minimum degree condition and for certain pseudo-random host graphs. See also Frieze and
Krivelevich~\cite{FriezeKrivelevich2002} and Chen, Chen, Im and Wang~\cite{ChenChenImWang2026} for
$(n,d,\lambda)$-graphs, and Wang~\cite{Wang2026} for balanced bipartite host graphs. For sparser host
graphs, Diskin and Geisler~\cite{DiskinGeisler2024} showed that if $\Gamma$ is $d$-regular with
$d\to\infty$ and $d\leq n^{1-a}$ for some fixed $a>0$, and satisfies suitable edge-expansion
conditions, then a.a.s.\ $\tau_k(\Gamma)=\tau_{k\text{-conn}}(\Gamma)$ for every fixed $k$. Closer to
our setting, Johansson~\cite{Johansson2020weighted} considered the inhomogeneous model in which each
edge $ij$ is present independently with probability $P_{ij}$, and proved hitting-time results for
edge-disjoint Hamilton cycles when $P$ satisfies a spectral gap condition and its minimum row sum is
of the same order as its average row sum.

Our setting differs from those considered above: the host graph itself is sampled from a graphon, and its edges are then revealed in a uniformly random order. We establish the coincidence of the connectivity and minimum-degree-one hitting times without imposing any quantitative expansion assumption. Instead, we require only that the graphon be connected and that the lower tail of its degree distribution satisfy the condition necessary for asymptotic connectivity. In particular, our result applies even when the graphon degree function is not bounded away from zero.

\subsection{Our contributions}
The first contribution of this paper is conceptual: we show that the lower-tail hypotheses of \Cref{thm:hladky} admit equivalent formulations in terms of vague convergence of rescaled measures and truncated Laplace transforms of the graphon degree distribution. Although these equivalences are relatively elementary, the resulting change of perspective is central to our approach: it allows us to address questions concerning low-degree counts, connectivity, and hitting times using Laplace-transform techniques, including the Laplace principle when suitable lower-tail asymptotics are available. To state this result precisely, we introduce the following definitions. For a uniform random variable $X\in [0,1]$ and a measurable set $B\subseteq[0,\infty)$, define the measure
\[
\nu_n(B):=n\cdot\PP(n d_W(X)\in B).
\]
Also, for fixed $K>0$, define the \emph{truncated Laplace transform} of the degree function as 
\begin{equation}\label{eq:truncLaplace}
\mathcal L_n(K):=n\cdot\EE\big[e^{-nd_W(X)}\Indic_{\{n d_W(X)\le K\}}\big].
\end{equation}
Recall that a sequence of measures $(\mu_n)_{n\in\mathbb N}$ converges~\emph{vaguely} to some measure $\mu $ if $\int fd\mu_n\to \int fd\mu$ for every continuous function $f:[0,+\infty)\to \mathbb R$ with compact support. Our first result is the following. 
\begin{theorem}\label{thm:measure-lowdeg}For a graphon $W$ and a fixed constant $c\geq 0$, the following statements are equivalent. 
\begin{enumerate}[label=(\alph{enumi})]
    \item\label{thm:measure-lowdeg:itm1} $\displaystyle\lim_{t\to 0}\frac{g_W(t)}{t}=c$.
    \item\label{thm:measure-lowdeg:itm2} $\nu_{n}\underset{n}{\longrightarrow} c\lambda$ in the vague topology.
    \item\label{thm:measure-lowdeg:itm3} For every fixed $K>0$, $\mathcal L_n(K)\underset{n}{\longrightarrow} c(1-e^{-K})$.
\end{enumerate}\end{theorem}
Moreover, we show that asymptotic upper and lower bounds on $g_W(u/n)/(u/n)$ as $n\to\infty$ are equivalent to corresponding bounds on $\nu_n([0,u])$ and $\mathcal L_n(K)$, for every fixed $u,K>0$ (see \Cref{prop:laplace-upper,prop:laplace-lower}). Using this reformulation, 
we obtain a simplified proof of all four parts of \Cref{thm:hladky}, together with corresponding results for rescaled graphons (see \Cref{thm:triangular}). Moreover, it allows us to clarify the
regime where $\lim_{t\to 0}\frac{g_W(t)}{t}=c$ for some $c\in (0,+\infty)$, which is not covered by \Cref{thm:hladky} and was explicitly asked in~\cite[Section~4.2]{HV}. Hladk\'y and Viswanathan~\cite{HV} showed that 
$\lim_{t\to 0}\frac{g_W(t)}t=c$  alone does not determine the
limiting distribution of the minimum degree.  Indeed, while the lower tail $g_W$ describes the abundance of low-degree vertex types, it provides no information about the overlap of their neighborhoods. Two graphons with the same lower-tail degree asymptotics may therefore exhibit different limiting distributions of their low-degree vertices. 

We identify an additional second-order decorrelation condition under which the counts of vertices of each fixed degree converge jointly to independent Poisson random variables, each with mean $c$, and consequently the minimum degree converges to a geometric law. Our second-order decorrelation condition controls precisely this overlap: it requires the aggregate codegree of pairs of types with degrees at most $K/n$, suitably rescaled, to vanish as $n\to\infty$, for every fixed $K>0$. Perhaps surprisingly, this requirement admits two equivalent analytic formulations: the vague convergence to zero of a measure on pairs of rescaled degrees, and the vanishing of its associated truncated Laplace transforms. Moreover, we show that a sufficient condition for this second-order decorrelation is that the normalized overlap between the neighborhoods of low-degree types vanishes uniformly as their degrees approach zero. Together with the first-order condition on $g_W$, this additional control yields independent Poisson limits for the counts of vertices of each fixed degree, and consequently a geometric limiting distribution for the minimum degree. 

To formulate this sufficient condition, we need the following definitions. Recall that for a measurable function $f:A\to\mathbb R$ on a measurable set $A$ of positive measure, we write
$\operatorname*{ess\,inf}_{A}f:=\sup\{a\in\mathbb R:f\ge a\text{ almost everywhere on }A\}$ and
$\operatorname*{ess\,sup}_{A}f:=\inf\{b\in\mathbb R:f\le b\text{ almost everywhere on }A\}$. Given a graphon $W:[0,1]^2\to [0,1]$,  the \emph{codegree kernel} defined by $W$ is the function
\begin{equation}\label{def:codegree}
q_W(x,y):=\int_0^1 W(x,z)W(y,z)\,dz\quad\text{for $x,y\in [0,1]$},
\end{equation}
which measures the overlap between the neighbors of vertices of types $x$ and $y$ in $G(n,W)$. For $t>0$, let~$B_W(t)$ be the set of pairs $(x,y)$ such that $x,y\in S_W(t)$ and $d_W(x)+d_W(y)>0$, and let $\eta:[0,1]\to \mathbb R$ be the function defined by 
\begin{equation}\label{def:eta}
\eta(t):=
\operatorname*{ess\,sup}_{\substack{(x,y)\in B_W(t)}}
\frac{q_W(x,y)}{d_W(x)+d_W(y)}
\end{equation}
if $B_W(t)$ has positive measure, and set $\eta(t):=0$ otherwise. Pairs with $d_W(x)+d_W(y)=0$ are omitted from the quotient, as for such pairs $q_W(x,y)=0$ and the inequalities based on \eqref{def:eta} will remain valid. We can now state our second result, which gives a sufficient condition for the convergence in distribution of $\delta(G(n,W))$.
\begin{theorem}\label{thm:distribution}Let $c\in (0,+\infty)$ and let $W$ be a graphon such that $\lim_{t\to 0}\frac{g_W(t)}{t}=c$ and $\lim_{t\to 0}\eta(t)=0$. Then, $\delta(G(n,W))$ converges in distribution to a geometric random variable on $\{0,1,2,\ldots\}$ with parameter $1-e^{-c}$.\end{theorem}
\Cref{thm:distribution} is  a corollary of a more general result (\Cref{thm:mixed-factorial})
about the joint distribution of the low-degree vertices in $G(n,W)$. Writing $N^{(r)}$ for the number
of vertices of $G(n,W)$ of degree exactly $r$, \Cref{thm:mixed-factorial}
shows that whenever $\lim_{t\to 0}\frac{g_W(t)}{t}=c$ and the
second-order condition is satisfied, the mixed
factorial moments of $(N^{(0)},\dots,N^{(\ell-1)})$ converge, so that $(N^{(0)},\dots,N^{(\ell-1)})$ converges jointly to independent
Poisson random variables of mean $c$.  In particular,  under these hypotheses, the number of isolated vertices in $G(n,W)$ converges to a Poisson distribution of mean $c$. \Cref{thm:distribution} thus
follows from the identity
$\{\delta(G(n,W)) \ge r\} = \{N^{(0)} = \ldots = N^{(r-1)} = 0\}$, which implies $\mathbb{P}(\delta(G(n,W)) \ge r) \to e^{-cr}$.

We now turn to the main focus of this paper: the random subgraph process on a $W$-random graph, which we call the \emph{$W$-random graph process} and define as follows. Let $W$ be a graphon, let $G\sim G(n,W)$, and let~$Y:=e(G)$ denote the number of edges in $G$. Let $e_1,\ldots,e_Y$ be a uniform random ordering of the edges of $G$ and, for $k\ge0$, let $G_k$ be the graph consisting of the first $\min\{k,Y\}$ revealed edges of $G$. For a fixed integer~$\ell\ge1$, the $\ell$-minimum-degree hitting time is 
\[
\tau_\ell(W):=\min\{k:\delta(G_k)\ge\ell\},\]
with the convention that $\tau_\ell(W)=+\infty$ if no $G_k$ has minimum degree at least $\ell$. Similarly, we let
\[\tau_{\mathrm{conn}}(W):=\min\{k:G_k~\text{ is connected}\}\]
denote the hitting time for connectivity, with the convention that $\tau_{\mathrm{conn}}(W)=+\infty$ whenever no $G_k$ is connected. To state our next result, we need the notion of graphon connectivity. We say that $W$ is \emph{connected} if there is no measurable set $A\subseteq[0,1]$ with $0<\lambda(A)<1$ such that $W=0$ almost everywhere on $A\times([0,1]\setminus A)$. 
\begin{theorem}\label{theorem:hitting-time}Let $W$ be a connected graphon and suppose that $\lim_{t\to 0}\frac{g_W(t)}{t}=0$. Then, 
\[
\lim_{n\to\infty}\PP(\tau_{\mathrm{conn}}(W)=\tau_1(W)<+\infty)=1.
\] \end{theorem}
We observe that both hypotheses in \Cref{theorem:hitting-time} are needed. Indeed, $\tau_1(W)<+\infty$ if and only if~$G(n,W)$ has no isolated vertices, and $\tau_{\mathrm{conn}}(W)<+\infty$ if and only if $G(n,W)$ is connected. If $W$ is not connected and $\lim_{t\to 0}\frac{g_W(t)}{t}=0$, then $G(n,W)$ is a.a.s.\ disconnected and has no isolated vertices, in which case $\tau_{\mathrm{conn}}(W)=+\infty$ but $\tau_1(W)<+\infty$. On the other hand, if $\liminf_{t\to 0}\frac{g_W(t)}{t}>0$, then $G(n,W)$ has an isolated vertex with positive probability and so $\tau_1(W)=\tau_{\mathrm{conn}}(W)=+\infty$. 

We remark that the process considered here is not
identical to the \emph{$W$-edge-incremental process} proposed by Hladk\'y and Viswanathan~\cite{HV}, which is defined as follows. Given a graphon $W$, let $x_1,\ldots, x_n$ be the types produced in the vertex-generation stage. The $W$-edge incremental process is a sequence of graphs $(G^t)_{t=0}^{N}$ on the vertex set $[n]$, where $N:=\binom{n}{2}$ and $G^0:=\emptyset$ is the empty graph. For $t\ge 1$, say a pair $ij\not\in E(G^{t-1})$ is \textit{visible} if $W(x_i,x_j)>0$. If there are visible pairs, then $G^t$ is formed from $G^{t-1}$ by adding a random visible pair $ij$ with probability  $W(x_i,x_j)/(\sum_{ab\in \binom{[n]}{2}\setminus E(G^{t-1})}W(x_a,x_b))$. Otherwise, $G^t$ is formed from $G^{t-1}$ by adding a uniform non-edge. In particular, this process always ends with the complete graph $K_n$, so its hitting times are always finite, while the $W$-random graph process ends with $G(n,W)$. Hladk\'y and Viswanathan conjectured that, for every connected graphon $W$ and every fixed $K\ge 1$, the hitting time for $K$-connectivity and having minimum degree at least $K$ coincide in the $W$-edge-incremental process. Thus \Cref{theorem:hitting-time} is an analogue for this conjecture in the $W$-random graph process for $K=1$.

Finally, we also study the location and sharpness of the common hitting time $\tau_{\mathrm{conn}}(W)=\tau_{1}(W)$. For a connected graphon $W$, we show that there is a deterministic time $k=\Theta(n\log n)$ by which every fixed positive minimum degree is attained a.a.s.\ if and only if $d_*:=\operatorname*{ess\,inf}d_W>0$ (see Proposition~\ref{prop:nlogn-essential-inf}). Moreover, we show that if $d_*>0$, then a.a.s.\ $\tau_1(W)=\tau_{\mathrm{conn}}(W)=(1+o(1))\frac{\|W\|_1}{2d_*}n\log n$, which coincides with the case of $G(n,p)$ for any fixed $p\in (0,1]$ (see \Cref{prop:precise hitting time}). That is,  $G(n,W)$ behaves like the random graph process when $d_*>0$. On the other hand, we show that a positive essential
infimum is not necessary for a multiplicatively sharp transition.  That is, we give a
simple analytic condition on the Laplace transform of the graphon degree distribution under which there are deterministic
 times $k_-< k_+$ with $k_-/k_+\rightarrow 1$ such that a.a.s.\ $k_-\le \tau_1(W)\le k_+$, and exhibit
connected graphons with $\operatorname*{ess\,inf}d_W=0$ whose
connectivity transition is nevertheless multiplicatively sharp (see \Cref{prop:sharp-window}). 

A central feature of our approach is the systematic use of Laplace transforms of the graphon degree distribution to study connectivity and hitting times. By expressing expected counts of low-degree vertices in terms of Laplace-type integrals, we reduce several probabilistic questions to the analysis of these integrals through the lower tail $g_W$. When the lower tail has a suitable asymptotic form, the Laplace transform allows us to estimate these integrals by minimizing the corresponding exponent. In particular, this provides a way to determine the location and sharpness of connectivity transitions, including those occurring at scales larger than $n\log n$. Together with our comparison between the edge-revelation and thinning processes, this yields an analytic approach to hitting-time questions in inhomogeneous random graphs.

\subsection{Organization}
In \Cref{sec:analytic} we prove \Cref{thm:measure-lowdeg} and we use the codegree kernel to define a second-order decorrelation condition. \Cref{sec:static} revisits the inhomogeneous random graph model, gives shorter proofs of the main minimum-degree and connectivity results of~\cite{HV}, proves the Poisson limit under the second-order decorrelation condition and finishes with the proof of \Cref{thm:distribution}. \Cref{sec:dynamic} collects the results concerning the~$W$-random graph process, including estimates for exact-degree counts and minimum-degree hitting times, as well as our connectivity hitting-time theorem (\Cref{theorem:hitting-time}). We conclude with open problems in \Cref{sec:discussion}.

\section{Low graphon degrees: tails, Laplace transforms, and vague limits}\label[section]{sec:analytic}

In this section, we establish the relations between the lower tail of the graphon degree distribution, the rescaled family of measures $\nu_n$, and the
truncated Laplace transform. These results give a measure-theoretic interpretation of
the hypotheses in Theorem~3 of Hladk\'y and Viswanathan~\cite{HV} (see \Cref{thm:hladky}) and will be key when analyzing the $W$-random graph process.  Given a graphon $W$, recall that $\nu_n$ is defined by 
\[\nu_n(B)=n\cdot\mathbb P(nd_W(X)\in B),\]
where $X\sim U([0,1])$  and $B\subseteq [0,+\infty)$ is a measurable set. Note that for every $u\ge 0$, the definition of $\nu_n$ gives
\[\nu_n([0,u])=n g_W(u/n).\]
Also, for $K>0$, the truncated Laplace transform defined in~\eqref{eq:truncLaplace} satisfies
\[\mathcal L_n(K)=n\int_{S_W(K/n)}e^{-nd_W(x)}\,dx
      =\int_{[0,K]}e^{-u}\,\nu_n(du).
\]
Thus, $\mathcal L_n(K)$ is a rescaled truncated Laplace--Stieltjes transform of the distribution of $d_W(X)$. Moreover, Stieltjes integration by parts gives
\begin{equation}\label{eq:truncated-laplace-ibp}
\mathcal L_n(K)=e^{-K}\nu_n([0,K])+\int_0^K e^{-u}\nu_n([0,u])\,du.
\end{equation}
With these observations at hand, we can prove an upper-bound equivalence as in \Cref{thm:measure-lowdeg}.

\begin{proposition}[Upper one-sided bounds]\label[proposition]{prop:laplace-upper}
For $C\in[0,\infty)$ fixed, the following are equivalent.
\begin{enumerate}[label=(\roman*)]
\item\label{itm1:prop:up}
$\displaystyle\limsup_{t\to 0}\frac{g_W(t)}t\le C$.

\item\label{itm2:prop:up} For every fixed $u>0$,
$\displaystyle\limsup_{n\to\infty}\nu_n([0,u])\le Cu$.

\item\label{itm3:prop:up} For every fixed $K>0$,
$\displaystyle\limsup_{n\to\infty}\mathcal L_n(K)\le C(1-e^{-K})$.
\end{enumerate}
Moreover, under these conditions, $(\nu_n)_{n\in\mathbb N}$ is relatively compact in the vague topology, meaning that every subsequence of $(\nu_n)_{n\in\mathbb N}$ has a subsequence that converges vaguely to a locally finite measure $\nu$, which satisfies~$\nu([0,u])\le Cu$ for every $u>0$.
\end{proposition}

\begin{proof} We assume first~\cref{itm1:prop:up}. For $u>0$, note that we have
$\nu_n([0,u])=u\,g_W(u/n)/(u/n)$,
and therefore $\limsup_n \nu_n([0,u])\le Cu$. This proves~\cref{itm2:prop:up}.

Assume now~\cref{itm2:prop:up} and fix $K>0$. Since $u\mapsto\nu_n([0,u])$ is nondecreasing and $\limsup_n \nu_n([0,K])\le CK$, the functions $u\mapsto\nu_n([0,u])$ are uniformly bounded on $[0,K]$ for all sufficiently large $n$. Then, by the reverse Fatou inequality applied to~\eqref{eq:truncated-laplace-ibp}, we have
\begin{align*}
\limsup_{n\to\infty}\mathcal L_n(K)
&\le
e^{-K}\limsup_{n\to\infty}\nu_n([0,K])
+
\int_0^K e^{-u}\limsup_{n\to\infty}\nu_n([0,u])\,du\\
&\le
C\Big(Ke^{-K}+\int_0^Kue^{-u}\,du\Big) \\
&=
C(1-e^{-K}).
\end{align*}
This proves~\cref{itm3:prop:up}. Finally, let us assume~\cref{itm3:prop:up}. Let $K>0$ be arbitrary. Note that, as $e^{-u}\ge e^{-K}$ for all $u\in [0,K]$, we have $\mathcal L_n(K)\ge e^{-K}\nu_n([0,K])$. Therefore,
\[
\limsup_{n\to\infty}\frac{\nu_n([0,K])}{K}
\le
C\frac{e^K-1}{K}.
\]
For $t>0$, choose $n$ so that $\frac{K}{n+1}<t\le \frac{K}n$. Then, by monotonicity of $g_W$, we have
\[\frac{g_W(t)}{t}\le
\frac{n+1}{K}g_W(K/n)
=
\frac{n+1}{n}\frac{\nu_n([0,K])}{K}.
\]
As $n\to \infty$ when $t\to 0$, it follows that
\[
\limsup_{t\to 0}\frac{g_W(t)}{t}
\le
C\frac{e^K-1}{K}.
\]
Hence letting $K\to0$ proves~\cref{itm1:prop:up}.

It remains to prove the second assertion. First, to see that $(\nu_n)_{n\in\mathbb N}$ is relatively compact, it suffices to show that $\sup_{n\in\mathbb N}\nu_n([0,K])<\infty$ for every $K>0$. Indeed, for a fixed
$K>0$, by \cref{itm2:prop:up} there is some~$n_0$ such that $\sup_{n\ge n_0}\nu_n([0,K])<\infty$, and since each
$\nu_1,\ldots,\nu_{n_0-1}$ also has finite mass on $[0,K]$, it follows that $\sup_{n\in\mathbb N}\nu_n([0,K])<\infty$. By the standard compactness criterion for locally finite measures, uniform boundedness of the masses on compact sets implies that $(\nu_n)_{n\in\mathbb N}$ is relatively compact in the vague topology.

Now suppose $(\nu_{n_j})_{j\in\mathbb N}$ converges vaguely to some measure $\nu$. Fix $u>0$ and $\delta>0$, and choose a continuous compactly supported function $\phi:[0,+\infty)\to \mathbb R$ such that $0\le\phi\le1$, $\phi\equiv1$ on $[0,u]$, and $\phi\equiv 0$ on $[u+\delta,\infty)$. Then,
\[
\nu([0,u])
\le
\int\phi\,d\nu
=
\lim_{j\to\infty}\int\phi\,d\nu_{n_j}
\le
\limsup_{j\to\infty}\nu_{n_j}([0,u+\delta])
\le
C(u+\delta).
\]
Letting $\delta\to0$ thus finishes the proof.
\end{proof}

The analogous statement holds for the lower bounds.

\begin{proposition}[Lower one-sided bounds]\label[proposition]{prop:laplace-lower}
For $c\ge 0$ fixed, the following are equivalent.
\begin{enumerate}[label=(\roman*)]
\item\label{itm1:prop:low}
$\displaystyle\liminf_{t\to 0}\frac{g_W(t)}{t}\ge c$.
\item\label{itm2:prop:low} For every fixed $u>0$,
$\displaystyle\liminf_{n\to\infty}\nu_n([0,u])\ge cu$.
\item\label{itm3:prop:low} For every fixed $K>0$, 
$\displaystyle\liminf_{n\to\infty}\mathcal L_n(K)\ge c(1-e^{-K})$.
\end{enumerate}
Moreover, if $(\nu_{n_j})_{j\in\mathbb N}$ converges vaguely to a locally finite measure $\nu$, then
$\nu([0,u])\ge cu$
for every $u>0$.
\end{proposition}

\begin{proof}
For $c=0$, all assertions follow immediately from nonnegativity. Hence, assume $c>0$. 
The proof follows the same argument as the proof of \Cref{prop:laplace-upper}, so we only highlight the minor differences. To deduce~\cref{itm1:prop:low} from~\cref{itm2:prop:low}, choose $n$ such that $K/n\le t<K/(n-1)$ so
\[
\frac{g_W(t)}{t}
\ge
\frac{n-1}{K}g_W(K/n)
=
\frac{n-1}{n}\frac{\nu_n([0,K])}{K},
\]
and conclude by adapting the argument of \Cref{prop:laplace-upper}.

For the moreover part, suppose $\nu_{n_j}\to\nu$ vaguely. Fix $u>0$ and $\delta\in(0,u)$, and choose a continuous function $\phi:[0,+\infty)\to\mathbb R$ supported on $[0,u]$ such that $0\le\phi\le1$, $\phi\equiv1$ on $[0,u-\delta]$, and $\phi\equiv0$ on $[u,\infty)$. Then,
$\nu([0,u])\ge\int\phi\,d\nu=\lim_{j\to\infty}\int\phi\,d\nu_{n_j}
\ge\liminf_{j\to\infty}\nu_{n_j}([0,u-\delta])\ge c(u-\delta)$. As $\delta>0$ was arbitrary, we may let $\delta\to 0$ to finish the proof.
\end{proof}

Combining these two propositions, we may give the proof of \Cref{thm:measure-lowdeg}.

\begin{proof}[Proof of \Cref{thm:measure-lowdeg}]
Assume first that $\lim_{t\to 0}\frac{g_W(t)}t=c$. Then, by \Cref{prop:laplace-upper}, the sequence $(\nu_n)_{n\in\mathbb N}$ is relatively compact for the vague topology. Moreover, if $(\nu_{n_j})_{j\in\mathbb N}$ converges to some $\nu$, \Cref{prop:laplace-upper,prop:laplace-lower} imply that $\nu([0,u])=cu$ for every $u>0$. Hence, $\nu=c\lambda$ for every accumulation point of $(\nu_n)_{n\in\mathbb N}$, and therefore $\nu_n\to c\lambda$ in the vague topology.

Conversely, if $\nu_n\to c\lambda$ vaguely, then
$\nu_n([0,u])\to cu$ for every fixed $u>0$, since
$c\lambda$ assigns zero measure to the boundary of $[0,u]$.
Moreover, for every fixed $K>0$, the functions
$u\mapsto\nu_n([0,u])$ are uniformly bounded on $[0,K]$
for all sufficiently large $n$, because
$\nu_n([0,u])\le\nu_n([0,K])\to cK$.
Hence, \eqref{eq:truncated-laplace-ibp} and the dominated
convergence theorem give
\[
\lim_{n\to\infty}\mathcal L_n(K)=c\Big(Ke^{-K}+\int_0^Kue^{-u}\,du\Big)=c(1-e^{-K})
\]
for every fixed $K>0$. Finally, if~\cref{thm:measure-lowdeg:itm3} holds, its upper and lower bounds satisfy~\cref{itm3:prop:up,itm3:prop:low} 
of \Cref{prop:laplace-upper,prop:laplace-lower}
with constant $c$ and thus
\[\limsup_{t\to 0}\frac{g_W(t)}t\le c\le\liminf_{t\to 0}\frac{g_W(t)}t,\]
thus finishing the proof.
\end{proof}

\subsection{The codegree kernel and the second-order decorrelation condition}\label[section]{sec:second-order-vague}

Though the measure $\nu_n$ records the expected number of vertex types
whose graphon degrees are of order $O(1/n)$, it does not describe the overlap
between the neighborhoods of two such types, which is typically needed for a second-moment argument. To capture this information, we will use the codegree kernel $q_W$ defined in \cref{def:codegree}. For $t\ge0$, set
\[
J(t):=\int_{S_W(t)^2}q_W(x,y)\,dx\,dy,
\]
and, as with $\nu_n$, we define a measure $\kappa_n$ on $[0,\infty)\times [0,\infty)$ by
\[
\kappa_n(B):=
n^3\int_{[0,1]^2}
\Indic_{\{(nd_W(x),nd_W(y))\in B\}}q_W(x,y)\,dx\,dy.
\]
In particular, we have that $\kappa_n([0,K]^2)=n^3J(K/n)$ for every $K>0$.  We say that the low-degree types satisfy the \emph{second-order  decorrelation} condition if $\kappa_n\to 0$ in the vague topology. Likewise for the measure $\nu_n$, we will use a truncated Laplace transform for the measure $\kappa_n$, defined as
\[\mathcal L_n^{(2)}(K):=
\int_{[0,K]^2}e^{-(u+v)}\,d\kappa_n
=n^3\int_{S_W(K/n)^2}
 e^{-n(d_W(x)+d_W(y))}q_W(x,y)\,dx\,dy.
\]

\begin{proposition}\label[proposition]{prop:second-order-vague}
The following are equivalent:
\begin{enumerate}[label=(\roman*)]
\item\label{itm1:prop:second-order} $\kappa_n\xrightarrow[n]{\mathrm{vague}}0$ on $[0,\infty)^2$.
\item\label{itm2:prop:second-order} For every fixed $K>0$, $n^3J(K/n)\xrightarrow[n]{}0$.
\item\label{itm3:prop:second-order} $J(t)/t^3\to 0$ as $t\to0$.
\item\label{itm4:prop:second-order} For every fixed $K>0$, $\mathcal L_n^{(2)}(K)\xrightarrow[n]{}0$.
\end{enumerate}
\end{proposition}

\begin{proof}
Because the measures $\kappa_n$ are positive, vague convergence to zero is equivalent to the convergence to zero of their masses on every fixed compact square. Indeed, if $\kappa_n\to0$ vaguely, a continuous compactly supported function that is at least one on $[0,K]^2$ shows that $\kappa_n([0,K]^2)\to0$. Conversely, if these masses tend to zero for every fixed $K$, then every compactly supported continuous test function is dominated by a constant times the indicator of some $[0,K]^2$. Since
$\kappa_n([0,K]^2)=n^3J(K/n)$, 
this proves the equivalence of~\cref{itm1:prop:second-order} and~\cref{itm2:prop:second-order}.

On $[0,K]^2$ one has $e^{-2K}\le e^{-(u+v)}\le1$, and hence
\[
e^{-2K}\kappa_n([0,K]^2)
\le \mathcal L_n^{(2)}(K)
\le \kappa_n([0,K]^2).
\]
Thus~\cref{itm2:prop:second-order} and~\cref{itm4:prop:second-order} are equivalent.

If~\cref{itm3:prop:second-order} holds, then for every fixed $K>0$,
\[
n^3J(K/n)=K^3\frac{J(K/n)}{(K/n)^3}\xrightarrow[n\to\infty]{}0,
\]
so~\cref{itm2:prop:second-order} follows. Conversely, assume~(ii). Given $t>0$, choose $n$ so that $1/(n+1)<t\le1/n$. Since $J$ is nondecreasing,
\[
\frac{J(t)}{t^3}
\le (n+1)^3J(1/n)
=\Big(\frac{n+1}{n}\Big)^3 n^3J(1/n)
\xrightarrow[n\to\infty]{}0.
\]
This proves~\cref{itm3:prop:second-order}.
\end{proof}
Note that for two types $x$ and $y$, the quantities $d_W(x)$ and $d_W(y)$ measure
the probabilities of being adjacent to a uniformly chosen third
type, whereas
\[
q_W(x,y)=\int_0^1 W(x,z)W(y,z)\,dz
\]
is the probability that this third type is adjacent to both. In the
binomial model $G(n,p)$, the corresponding probabilities are $p$ and
$p^2$, respectively. Thus, in the low-degree regime, the
natural analogue of this independence is that the common-neighbor intensity
of two low-degree types is negligible compared with their individual
neighbor intensities, which motivates the condition $\eta(t)\to 0$ introduced in \eqref{def:eta}.

The next corollary connects this combinatorial notion of asymptotic
independence with the measure-theoretic condition introduced above. Under the
linear upper bound on the low-degree tail needed for connectivity, the uniform decay of the normalized
common-neighbor overlap implies the second-order  decorrelation condition.
\begin{corollary}\label[corollary]{cor:eta-implies-second-order}
Suppose $g_W(t)=O(t)$ as $t\to0$ and $\eta(t)\to0$ as $t\to0$. Then,
$\kappa_n\xrightarrow[n\to\infty]{\mathrm{vague}}0.$
\end{corollary}

\begin{proof}
For $x,y\in S_W(t)$, the definition of $\eta(t)$ gives $q_W(x,y)\le\eta(t)(d_W(x)+d_W(y))$ for $t$ sufficiently small. Therefore, using that $d_W(x)\le t$ for $x\in S_W(t)$, we have
\[
J(t)
\le2\eta(t)g_W(t)\int_{S_W(t)}d_W(x)\,dx
\le2\eta(t)t g_W(t)^2
=o(t^3).
\]
Hence the conclusion follows from \Cref{prop:second-order-vague}.
\end{proof}

\section{Minimum degree of $G(n,W)$ revisited}\label[section]{sec:static}
In this section, we apply the analytic framework of \Cref{sec:analytic} to study the minimum degree of $G(n,W)$. 
In particular, we recover all four parts of \Cref{thm:hladky} from direct moment estimates for the numbers of
vertices of small degree. The absence of low-degree vertices is controlled by a first-moment argument, while the presence of low-degree vertices is established by a second-moment
argument after restricting to a suitably chosen set of types of low graphon degree.
Although the ingredients behind these estimates are present in the
arguments of Hladk\'y and Viswanathan, they do not organize the proof of their theorem
in this way. Moreover, making this connection explicit will be particularly useful when analyzing the $W$-random graph process. 

\subsection{An alternative proof of Theorem~\ref{thm:hladky}}\label{ssec:alternativeProof}
For the rest of this section, suppose $W$ is a fixed graphon and let $G\sim G(n,W)$. In order to study the low-degree vertices in $G(n,W)$, we will consider, for $r\ge 0$,  the integral
\[I_r(C):=\int_0^1d_W(x)^re^{-Cd_W(x)}\,dx.\]
In particular, when $r = 0$, we define the \emph{non-truncated Laplace transform} by
\[
\mathcal{L}(C):= I_0(C) = \EE[e^{-Cd_W(X)}].
\]
We begin with a technical lemma that allows us to control the integrals $I_r(C)$ in the context of \Cref{thm:measure-lowdeg}.
\begin{lemma}\label[lemma]{lem:Ir-small-tail}Suppose that $\lim_{t\to 0}\frac{g_W(t)}t=c$ for some
$c\in[0,\infty)$. Then, for every fixed integer $r\ge0$,
\[
C^{r+1}I_r(C)\xrightarrow[C\to\infty]{} c\,r!.
\]
In particular, if $c=0$, then $I_r(C)=o(C^{-r-1})$, and
$\lim_{C\to \infty}C\cdot \mathcal L(C)=0$.
\end{lemma}
\begin{proof}
Since $g_W(0)=0$ and $g_W(1)=1$, Stieltjes integration by parts, with $\varphi(t)=t^re^{-Ct}$, gives
\[
I_r(C)
=e^{-C}-\int_0^1\varphi'(t)g_W(t)\,dt
=e^{-C}+\int_0^1g_W(t)(Ct^r-rt^{r-1})e^{-Ct}\,dt.
\]
After the change of variables $u=Ct$, we obtain
\[
C^{r+1}I_r(C)
=C^{r+1}e^{-C}
+\int_0^C Cg_W(u/C)(u^r-ru^{r-1})e^{-u}\,du.
\]

Since $g_W(t)/t\to c$ as $t\to0$ and $g_W(t)\le1$, there exists $\Lambda>0$ such that $g_W(t)\le\Lambda t$ for all $t\in[0,1]$. Hence, for $0\le u\le C$, $Cg_W(u/C)\le\Lambda u$.
Moreover, $Cg_W(u/C)\to cu$ for every fixed $u>0$. Since
\[
\big|Cg_W(u/C)(u^r-ru^{r-1})e^{-u}\big|
\le\Lambda(u^{r+1}-ru^r)e^{-u},
\]
and the right-hand side is integrable on $[0,\infty)$, the dominated convergence theorem, extending the integrand by zero for $u>C$, gives
\[
\lim_{C\to\infty}C^{r+1}I_r(C)
=c\int_0^\infty(u^{r+1}-ru^r)e^{-u}\,du
=c\big((r+1)!-r\cdot r!\big)=c\,r!.
\]
Here we also used $C^{r+1}e^{-C}\to0$, and the terms with a factor of $r$ are absent when $r=0$.
\end{proof}
As mentioned in the introduction, let $N^{(r)}:=|\{j\in[n]:d_G(j)=r\}|$ denote the random variable that counts the number of vertices in $G$ with degree exactly $r$.  For $\ell\ge 1$, we  write
\[
N^{(<\ell)}:=\sum_{r=0}^{\ell-1}N^{(r)}
=|\{j\in[n]:d_G(j)<\ell\}|,
\]
for the number of vertices with fewer than $\ell$ neighbors, and define $N^{(\le\ell)}$ likewise. In particular, note that~$N^{(0)}=N^{(<1)}=N^{(\le 0)}$ is the number of isolated vertices in $G$. The next result gives the exact first moment computation for $N^{(r)}$ and upper and lower bounds for $\mathbb E[N^{(r)}]$ in terms of the integrals $I_r(C)$.

\begin{lemma}\label[lemma]{lem:first-moment-static}
For every fixed $r\ge0$,
\[
\EE[N^{(r)}]
=n\binom{n-1}{r}
\int_0^1 d_W(x)^r\Big(1-d_W(x)\Big)^{n-1-r}\,dx.
\]
Moreover, if $\delta=\delta(n)\in(0,1)$ and $C_{r,n}:=n-1-r$, then
\[
n\binom{n-1}{r}
\Big(I_r\Big(\tfrac{C_{r,n}}{1-\delta}\Big)
-e^{-C_{r,n}\delta/(1-\delta)}\Big)
\le \EE[N^{(r)}]
\le n\binom{n-1}{r}I_r(C_{r,n}).
\]
\end{lemma}

\begin{proof}Let $X_1,\ldots, X_n\in [0,1]$ be the types of $G(n,W)$. Note that, conditional on $X_j=x$, the degree of the vertex $j$ in $G$ has distribution $\operatorname{Bin}(n-1,d_W(x))$. Therefore, 
\[\mathbb E[N^{(r)}]=n\int_0^1\mathbb P(\operatorname{Bin}(n-1,d_W(x))=r)\,dx=n\binom{n-1}{r}
\int_0^1 d_W(x)^r\Big(1-d_W(x)\Big)^{n-1-r}\,dx.\]
The upper bound follows from the inequality $1-d_W(x)\le e^{-d_W(x)}$. For the lower bound, note first that~$\log(1-u)\ge-u/(1-u)\ge-u/(1-\delta)$ for every $u\in [0,\delta]$. Hence
\begin{align*}
\int_0^1 d_W(x)^r\Big(1-d_W(x)\Big)^{C_{r,n}}\,dx 
&\ge \int_{S_W(\delta)}d_W(x)^re^{-C_{r,n}d_W(x)/(1-\delta)}\,dx\\
&= I_r\Big(\tfrac{C_{r,n}}{1-\delta}\Big)-\int_{[0,1]\setminus S_W(\delta)}d_W(x)^re^{-C_{r,n}d_W(x)/(1-\delta)}\,dx\\
&\ge  I_r\Big(\tfrac{C_{r,n}}{1-\delta}\Big)-e^{-C_{r,n}\delta/(1-\delta)},\end{align*}
where we used that $d_W(x)\le1$ for all $x\in [0,1]$. 
\end{proof}
We next use a truncated second-moment argument to recover simultaneously the two statements concerning isolated vertices in \Cref{thm:hladky}. 
\begin{proposition}
\label[proposition]{prop:static-ac} \mbox

\begin{enumerate}[label=(\roman*)]
    \item If $\liminf_{t\to 0}\frac{g_W(t)}t>0$, then $\liminf_{n\to\infty}\PP(\delta(G(n,W))=0)>0$.
\item If $\lim_{t\to 0}\frac{g_W(t)}t=\infty$, then a.a.s.\ $\delta(G(n,W))=0$.
\end{enumerate}
\end{proposition}

\begin{proof}
Let $0<L<\liminf_{t\to 0}\frac{g_W(t)}t$ and fix $q>0$. As $g_W(q/n)\ge Lq/n$ for all sufficiently large $n$, we may choose a measurable set $S_n\subseteq S_W(q/n)$ with $\lambda(S_n)=Lq/n$. Let $V_n$ count the number of isolated vertices whose associated types lie in $S_n$. Since every $x\in S_n$ satisfies $d_W(x)\le q/n$, we have
\[\EE[V_n]\ge Lq\Big(1-\frac qn\Big)^{n-1}.\]
Write $V_n=\sum_{i=1}^n I_i$, where $I_i$ is the indicator that the vertex $i$ is isolated in $G(n,W)$ and $X_i\in S_n$. 
Hence,~$\EE[I_i]\ge\lambda(S_n)(1-\frac{q}{n})^{n-1}$.
For $i\ne j$, we only use the crude inequality
$I_iI_j\le\Indic_{\{X_i\in S_n,X_j\in S_n\}}$. Independence of the two types therefore gives
$\EE[I_iI_j]\le\lambda(S_n)^2$. Consequently,
\[
\EE[V_n^2]\le\EE[V_n]+n(n-1)\lambda(S_n)^2
\le\EE[V_n]+(Lq)^2.
\]
Combining the bounds $Lq(1-q/n)^{n-1}\le\EE[V_n]\le Lq$
with the second-moment inequality yields 
\[
\liminf_{n\to\infty}\PP(V_n>0)
\ge
\liminf_{n\to\infty}\frac{(\EE[V_n])^2}{\EE[V_n]+(Lq)^2}
\ge
\frac{e^{-2q}}{1+(Lq)^{-1}}>0,
\]
which proves the first statement. 

For the second statement, suppose that $g_W(t)/t\to\infty$ as $t\to 0$, in which case the constant $L$ may be chosen arbitrarily large. By choosing $q=q(L)$ so that $q\to 0$ and $Lq\to\infty$ as $L\to\infty$, the last displayed lower bound then tends to one, which proves the second statement.
\end{proof}
We now consider the regime where $g_W(t)/t\to0$ as $t\to 0$. 

\begin{proposition}\label[proposition]{prop:static-b}
 If $\lim_{t\to 0}\frac{g_W(t)}t=0$, then $\delta(G(n,W))\to\infty$ in probability.
\end{proposition}

\begin{proof}
Fix $\ell\ge 1$. For every $0\le r<\ell$, \Cref{lem:Ir-small-tail} applied with $C=n-1-r$ gives $I_r(n-1-r)=o(n^{-r-1})$. Then, by \Cref{lem:first-moment-static}, $\EE[N^{(r)}]=o(1)$ for every $r<\ell$, and hence $\EE[N^{(<\ell)}]=o(1)$. Since $\{\delta(G)<\ell\}=\{N^{(<\ell)}>0\}$, Markov's inequality proves the assertion.
\end{proof}
So far we have recovered parts (a)--(c) of \Cref{thm:hladky}, thus it only remains to prove part~(d). In order to handle this case, though, we require an additional ingredient. We need the following result by Hladk\'y and Viswanathan, which shows that we can remove a small exceptional set containing
the lowest-degree types while ensuring that every type outside this set
retains a fraction of its original graphon degree.  In order to state this result, we need the following definition. For a measurable $A\subseteq[0,1]$, write
\[
d_W(x,A):=\int_A W(x,y)\,dy
\]
for the graphon degree of a vertex of type $x$ into $A$. 
\begin{lemma}[{\cite[Lemma 7]{HV}}]\label[lemma]{lem:exceptional}
For every $n$ and every $c>0$, there exists a measurable set $S\subseteq[0,1]$ with
$S_W(2c/n)\subseteq S$ and $\lambda(S)\le3g_W(2c/n)$, such that every $x\notin S$
satisfies
\[
d_W(x,[0,1]\setminus S)\ge d_W(x)/4.
\]
\end{lemma}
We give a short proof of Lemma~\ref{lem:exceptional} in the appendix (see \Cref{app:static-simplifications}). We are now ready for part (d) of \Cref{thm:hladky}.

\begin{proposition}\label[proposition]{prop:static-d}
If $\limsup_{t\to 0}\frac{g_W(t)}t<\infty$, then
$\liminf_{n\to\infty}\PP\big(\delta(G(n,W))>\ell\big)>0$ for every fixed $\ell\ge 0$.
\end{proposition}

\begin{proof}
Let $\Psi>0$ and $t_0\in(0,1)$ such that $g_W(t)\le\Psi t$
for every $0<t\le t_0$, which exist since $\limsup_{t\to 0}\frac{g_W(t)}t<\infty$. Let $c>0$ be a constant to be chosen later, and let $S$ be the exceptional
set from \Cref{lem:exceptional}, which satisfies $S_W(2c/n)\subseteq S$ and $\lambda (S)\le 3g_W(2c/n)\le 6\Psi c/n$ for all $n$ sufficiently large.
Let $\mathcal S:=\{X_i\notin S$ for every $i\in[n]\}$.
Then,
\[\PP(\mathcal S)
=\big(1-\lambda(S)\big)^n
\ge\Big(1-\frac{6\Psi c}{n}\Big)^n,
\]
hence $\liminf_{n\to\infty}\PP(\mathcal S)>0$.

Letting $\overline{S}:= [0,1]\setminus S$, observe that by conditioning on $(X_1,\ldots, X_n)\in \mathcal S$, the
types remain independent and are uniform on $\overline S$. Moreover, if we 
further condition on $X_i=x\in\overline S$, the type of any other vertex is uniform
on $\overline S$, and therefore the conditional probability that the edge
joining these two vertices is present equals
$\widetilde d(x):=
\frac{d_W(x,\overline S)}{\lambda(\overline S)}$.
By \Cref{lem:exceptional}, we have
$d_W(x,\overline S)\ge d_W(x)/4$. Moreover, $x\notin S_W(2c/n)$, and hence~$d_W(x)>2c/n$.

Fix $0\le r\le\ell$. For all sufficiently large $n$, we have
$\lambda(\overline S)\ge1/2$ and $n-1-r\ge n/2$. Conditional on
$(X_1,\ldots, X_n)\in \mathcal S$ and $X_i=x$, the degree of the vertex $i$ has distribution
$\operatorname{Bin}(n-1,\widetilde d(x))$. Hence the probability that $d_G(i)=r$, conditioned on $(X_1,\ldots, X_n)\in \mathcal S$ and $X_i=x$, is
\[
\binom{n-1}{r}\widetilde d(x)^r
(1-\widetilde d(x))^{n-1-r}.
\]
Since $\lambda(\overline S)\ge1/2$ and
$d_W(x,\overline S)\le d_W(x)$, we have
$\widetilde d(x)\le2d_W(x)$. On the other hand,
$d_W(x,\overline S)\ge d_W(x)/4$ and
$\lambda(\overline S)\le1$, so
$\widetilde d(x)\ge d_W(x)/4$. Moreover, conditional on
$\mathcal S$, the density of $X_i$ on $\overline S$ is
$1/\lambda(\overline S)\le2$. Therefore, using
$1-u\le e^{-u}$ and $\binom{n-1}{r}\le n^r/r!$ and the fact that every $x\in\overline S$ satisfies
$d_W(x)>2c/n$, we obtain that
\[
\EE[N^{(r)}\mid\mathcal S]
\le
\frac{2^{r+1}n^{r+1}}{r!}
\int_{\{d_W>2c/n\}}
d_W(x)^r e^{-n d_W(x)/8}\,dx.
\]
We now bound the last integral. Let $m:=n/8$, $a:=2c/n$ and choose $c\ge4\ell$. Then, for every $0\le r\le\ell$, the function
$t^r e^{-m t}$ is decreasing on $[a,\infty)$ and for all sufficiently
large $n$, we also have $a<t_0$. Hence, as in the proof of Lemma~\ref{lem:Ir-small-tail}, Stieltjes integration by parts on $[a,t_0]$ and the bound
$g_W(t)\le\Psi t$ give
\[
\int_{\{a<d_W\le t_0\}}
d_W(x)^r e^{-m d_W(x)}\,dx
\le
e^{-m t_0}
+m\Psi \int_a^\infty t^{r+1}e^{-m t}\,dt.
\]
Also, since $d_W(x)\le1$, we have $\int_{d_W>t_0}d_W(x)^r e^{-m d_W(x)}\,dx\le e^{-m t_0}$, and moreover,
\[
n^{r+1}m\int_a^\infty t^{r+1}e^{-m t}\,dt
=
8^{r+1}\int_{c/4}^\infty u^{r+1}e^{-u}\,du=8^{r+1}\Gamma(r+2,c/4),
\]
where $\Gamma(s,z):=\int_z^\infty u^{s-1}e^{-u}\,du$
denotes the upper incomplete gamma function.
Thus, for each $0\le r\le\ell$, $\EE[N^{(r)}\mid\mathcal S]$ is bounded
by a constant depending only on $r$ and $\Psi$ times
$\Gamma(r+2,c/4)$,
plus a term that is exponentially small in $n$.
Since $\Gamma(r+2,c/4)\to0$ as $c\to\infty$ and there are only
finitely many values $0\le r\le\ell$, we may choose
$c$ large enough so that
$\limsup_{n\to\infty}
\EE[N^{(\le\ell)}\mid\mathcal S]<1/2.$ Then, Markov's inequality gives~$\PP\big(N^{(\le\ell)}=0\mid\mathcal S\big)>1/2$
for all sufficiently large $n$. Note that on this event,
$\delta(G(n,W))>\ell$. Therefore, 
\[\mathbb P(\delta(G(n,W))> \ell)\ge \mathbb P(\delta(G(n,W))> \ell\mid \mathcal S)\cdot \mathbb P(\mathcal S)>\frac{1}{2}\mathbb P(\mathcal S)\]
for every $n$ sufficiently large, which finishes the proof as $\liminf_{n\to \infty}\mathbb P(\mathcal S)>0$.\end{proof}

\subsection{The limiting distribution of the minimum degree}
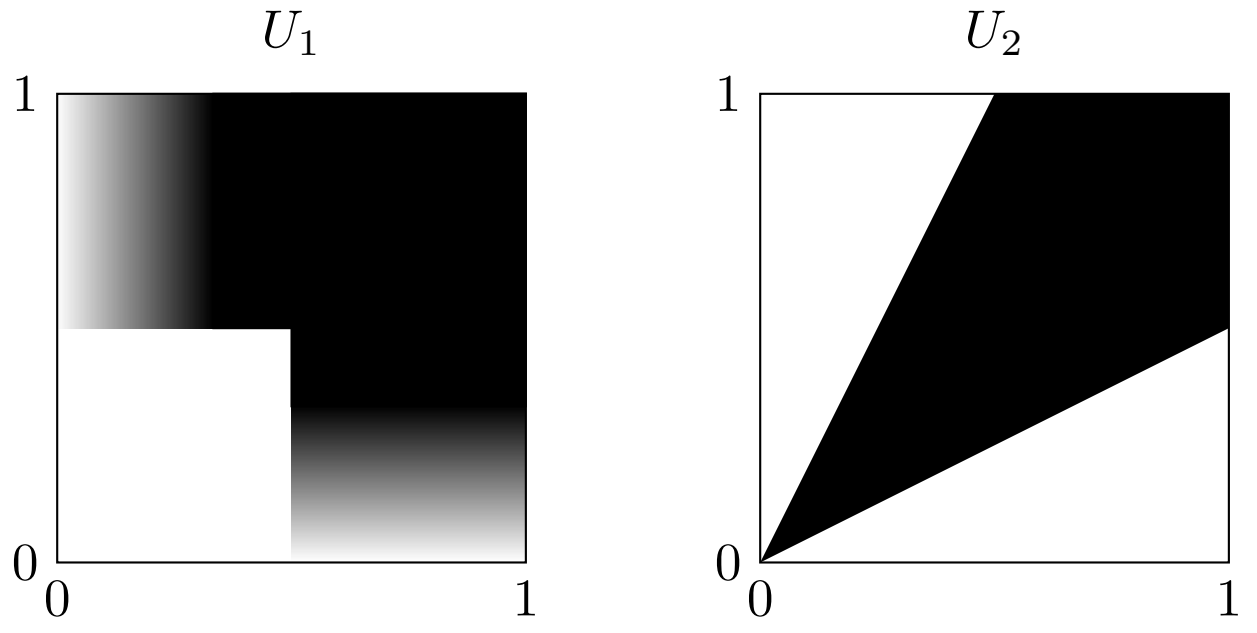
\begin{figure}[h!]
\centering
\begin{tikzpicture}[x=3.1cm,y=3.1cm]
  \begin{scope}
    \shade[left color=white,right color=black]
      (0,0.5) rectangle (0.333333,1);
    \fill[black] (0.333333,0.5) rectangle (0.5,1);

    \shade[bottom color=white,top color=black]
      (0.5,0) rectangle (1,0.333333);
    \fill[black] (0.5,0.333333) rectangle (1,0.5);

    \fill[black] (0.5,0.5) rectangle (1,1);

	\draw (0.5,0.333333) rectangle (1,1);
	\draw (0.333333,0.5) rectangle (1,1);
		
    \draw (0,0) rectangle (1,1);

    \node at (0.5,1.13) {$U_1$};
    \node[below] at (0,0) {$0$};
    \node[below] at (1,0) {$1$};
    \node[left] at (0,0) {$0$};
    \node[left] at (0,1) {$1$};
  \end{scope}

  \begin{scope}[xshift=4.65cm]
    \fill[black]
      (0,0) -- (1,0.5) -- (1,1) -- (0.5,1) -- cycle;
    \draw (0,0) rectangle (1,1);

    \node at (0.5,1.13) {$U_2$};
    \node[below] at (0,0) {$0$};
    \node[below] at (1,0) {$1$};
    \node[left] at (0,0) {$0$};
    \node[left] at (0,1) {$1$};
  \end{scope}
\end{tikzpicture}
\caption{Schematic representations of the graphons $U_1$ and $U_2$ from
Section~4.1 of~\cite{HV}. For $U_1$, white and black represent the values
$0$ and $1$, respectively, and the gray regions represent the linear
variation of the graphon before it reaches value $1$. The graphon $U_2$ is
the indicator of the black quadrilateral.}
\label{fig:U1-U2}
\end{figure}

\label[section]{sec:mixed-factorial}
We now turn our attention to the limiting distribution of $\delta(G(n,W))
$ under the assumption that $\lim_{t\to 0}\frac{g_W(t)}t=c$ for some $c\in(0,\infty)$. As we mentioned in the introduction, Hladk\'y and Viswanathan showed that this condition by itself is not sufficient to determine the limiting law of
the minimum degree. Specifically, they provided the counterexample described next. Let $U_1, U_2:[0,1]^2\to [0,1]$ be two graphons defined as follows. Let $U_1$ be zero on $[0,1/2)^2$,
equal to one on $[1/2,1]^2$, and given symmetrically on the two off-diagonal rectangles by
$U_1(x,y)=\min\{1,3\min\{x,y\}\}$. The graphon $U_2$ is the indicator of
the quadrilateral with vertices $(0,0)$, $(1,1/2)$, $(1,1)$, and
$(1/2,1)$, that is, $U_2(x,y):=\Indic_{\{x/2\le y\le2x\}}$ (see \Cref{fig:U1-U2}).

For $x\approx 0$, both graphons satisfy $d_{U_1}(x)=d_{U_2}(x)=3x/2$, and
hence $\lim_{t\to 0}\frac{g_{U_1}(t)}t=\lim_{t\to 0}\frac{g_{U_2}(t)}t=2/3$.
Nevertheless, Hladk\'y and Viswanathan show that
$\delta(G(n,U_1))$ converges to a geometric distribution with parameter
$1-e^{-2/3}$, whereas $\delta(G(n,U_2))$ converges to a geometric
distribution with parameter $1/2$. The difference is not visible in the degree distribution,
but rather in the way the neighborhoods of low-degree types overlap, which we measure through the codegree kernel $q_W$.
For $U_1$, if $x,y$ are sufficiently small, then
$q_{U_1}(x,y)=9xy/2$, so the codegree is much smaller
than the graphon degrees. By contrast, in $U_2$ the neighborhood
of a small type $x$ is the interval $[x/2,2x]$, and low-degree types of
comparable size have neighborhoods whose intersection has size comparable
to their individual degrees. Thus the two examples have the same
first-order lower tail but very different codegree kernels.

Recall that for a graphon $W$, $\nu_n$ and $\kappa_n$ are the measures defined by 
\[\nu_n(B)=n\int_{[0,1]}\Indic_{\{nd_W(x)\in B\}}\,dx\quad\text{and}\quad
\kappa_n(B)=
n^3\int_{[0,1]^2}
\Indic_{\{(nd_W(x),nd_W(y))\in B\}}q_W(x,y)\,dx\,dy.
\]
Before proving \Cref{thm:distribution}, we need the following technical lemma. 
\begin{lemma}\label[lemma]{lem:weighted-truncated-laplace}
Suppose that $\lim_{t\to 0}\frac{g_W(t)}t=c$ for some $c\in(0,\infty)$. Then, for every fixed $r\ge0$ and $K>0$,
\[
\frac{n}{r!}\int_{S_W(K/n)}(nd_W(x))^r e^{-nd_W(x)}\,dx
\xrightarrow[n\to\infty]{}
\frac{c}{r!}\int_0^K u^r e^{-u}\,du.
\]
\end{lemma}

\begin{proof}Recall that by \Cref{thm:measure-lowdeg} we have
$\nu_n\to c\,\lambda$ in the vague topology. Then, by the definition of $\nu_n$, we can write
\[
\frac{n}{r!}\int_{S_W(K/n)}(nd_W(x))^r e^{-nd_W(x)}\,dx
=
\frac1{r!}\int_{[0,K]}u^re^{-u}\,\nu_n(du).
\]
Since the limiting measure has no atoms, the vague convergence gives the result.
\end{proof}

We can now compute the mixed factorial moments of the low-degree count vector $(N^{(0)},\ldots,N^{(\ell-1)})$, from which \Cref{thm:distribution} will easily follow. We shall use the standard falling
factorial notation $(z)_a:=z(z-1)\cdots(z-a+1)$.

\begin{theorem}\label[theorem]{thm:mixed-factorial}Let $c\in(0,\infty)$ and suppose that $\nu_n\to c\,\lambda$ and $\kappa_n\to 0$, both in the vague topology. Fix $\ell\ge1$. Then, for every choice of fixed nonnegative integers $m_0,\ldots,m_{\ell-1}$,
\[
\EE\Big[\prod_{r=0}^{\ell-1}(N^{(r)})_{m_r}\Big]
\xrightarrow[n\to\infty]{}
c^{\,m_0+\cdots+m_{\ell-1}}.
\]
Consequently,
$(N^{(0)},\ldots,N^{(\ell-1)})
\xrightarrow[n\to\infty]{d}
(P_0,\ldots,P_{\ell-1})$,
where $P_0,\ldots,P_{\ell-1}$ are independent Poisson random variables with mean $c$.
\end{theorem}
\begin{proof}
Let $m:=\sum_{r=0}^{\ell-1}m_r$. If $m=0$, the assertion is immediate,
so assume $m\ge1$. Choose $r_1,\ldots,r_m\in\{0,\ldots,\ell-1\}$ so
that exactly $m_r$ of the $r_i$'s are equal to $r$. We call the vertices
in $[m]$ the distinguished vertices and those in
$\{m+1,\ldots,n\}$ the remaining vertices. By exchangeability, we have
\[
\EE\Big[\prod_{r=0}^{\ell-1}(N^{(r)})_{m_r}\Big]
=
(n)_m\PP\big(d_G(i)=r_i\text{ for every }i\in[m]\big).
\]
We estimate the probability on the right-hand side by conditioning on the types of
the distinguished vertices. We first compute the contribution from
configurations in which all of these types have graphon degree at most
$K/n$, for fixed $K>0$, and then show that the complementary contribution
vanishes as $K\to\infty$.

By \Cref{thm:measure-lowdeg}, $\nu_n\xrightarrow[n\to\infty]{\mathrm{vague}}c\,\lambda$ implies that
$\lim_{t\to 0}\frac{g_W(t)}t=c$. Fix $K>0$, set
$A_n:=S_W(K/n)$ and $a_n:=g_W(K/n)$. We condition on
$X_i=x_i$, $i\in[m]$, with $x:=(x_1,\ldots,x_m)\in A_n^m$, and write
\[
d_i:=d_W(x_i),\qquad
D(x):=\sum_{i=1}^m d_i,\qquad
Q(x):=\sum_{1\le i<j\le m}q_W(x_i,x_j).
\]
Since $g_W(t)/t\to c$, we have $a_n=O(1/n)$. Moreover,
\Cref{prop:second-order-vague} gives $J(K/n)=o(n^{-3})$, and hence,
when $m\ge2$, 
the probability that the
distinguished types all lie in $A_n$ and a fixed remaining vertex is adjacent to at least
two distinguished vertices is at most 
\begin{equation}\label{eqn:thm:mixed-factorial:Q}
\int_{A_n^m}Q(x)\,dx
=
\binom m2a_n^{m-2}J(K/n)
=
o(n^{-(m+1)}).
\end{equation}
Also, note that  $Q\equiv0$ for $m=1$.

We next describe the edges from the remaining vertices to the distinguished
ones. For $x\in A_n^m$, let $p_i(x)$ be the probability that a fixed
remaining vertex is adjacent to distinguished vertex $i$ and to none of the
other distinguished vertices, and let $p_0(x)$ be the probability that it
is adjacent to none of the distinguished vertices. For a fixed remaining
vertex, the probability of being adjacent to distinguished vertex $i$ is
$d_i$, while the probability of being adjacent to both $i$ and $j$ is
$q_W(x_i,x_j)$. Therefore, by the first two Bonferroni inequalities,
\[
0\le d_i-p_i(x)\le\sum_{j\ne i}q_W(x_i,x_j)
\quad\text{and}\quad
1-D(x)\le p_0(x)\le1-D(x)+Q(x).
\]
Let $R:=\sum_{i\in [m]} r_i$. For a fixed remaining vertex, the
probability of being adjacent to at least two distinguished vertices is at
most $Q(x)$, so the probability that this occurs for at least one of the
$n-m$ remaining vertices is at most $(n-m)Q(x)\le nQ(x)$. If no remaining vertex has two distinguished neighbors, the probability that, for every $i\in[m]$, exactly $r_i$
remaining vertices are adjacent to the distinguished vertex $i\in [m]$ is
\begin{equation}\label{eqn:multinom}
\frac{(n-m)_R}{\prod_i r_i!}
p_0(x)^{n-m-R}\prod_{i=1}^m p_i(x)^{r_i}.
\end{equation}
Let $\mathcal I$ denote the event that there is no edge between two
distinguished vertices. Conditional on the types of the distinguished
vertices, $\mathcal I$ is independent of all edges joining a distinguished
vertex to a remaining vertex. Conditional on the distinguished types and on $\mathcal I$, the joint degree probability differs from the last displayed multinomial expression by at most $nQ(x)$.

We now approximate the multinomial probability in~\eqref{eqn:multinom} by a product of Poisson
probabilities.
Since $d_i\le K/n$, we have $D(x)\le mK/n$. As $0\le np_i(x)\le nd_i\le K$ and the $r_i$'s are fixed, the
terms $(np_i(x))^{r_i}$ and~$(nd_i)^{r_i}$ are bounded by a constant $C_1=C_1(K,\ell)>0$.
Telescoping the difference
of the two products gives
\begin{align*}
\Big|
\prod_{i=1}^m (np_i(x))^{r_i}
-
\prod_{i=1}^m (nd_i)^{r_i}
\Big| 
& =
\Big|\sum_{j=1}^m
\Big(\prod_{i<j} (np_i(x))^{r_i}\Big)
\big((np_j(x))^{r_j}-(nd_j)^{r_j}\big)
\Big(\prod_{i>j} (nd_i)^{r_i}\Big)
\Big| \\
&  \le C_1^{m-1}\sum_{j=1}^m\Big|\big((np_j(x))^{r_j}-(nd_j)^{r_j}\big)\Big|. 
\end{align*}
Next, because $u\mapsto u^{r_i}$ is Lipschitz on $[0,K]$,
\[
\big|(np_j(x))^{r_j}-(nd_j)^{r_j}\big|
\le C_2n|p_j(x)-d_j|
\le C_2n\sum_{i\ne j}q_W(x_i,x_j),
\]
for some constant $C_2=C_2(K,\ell)$. This also holds when $r_j=0$, since
then both powers equal $1$.
It follows that, 
\begin{equation}\label{eq:midStep}
\Big|
\prod_{i=1}^m (np_i(x))^{r_i}
-
\prod_{i=1}^m (nd_i)^{r_i}
\Big|
\le 
2C_1^{m-1}C_2nQ(x).
\end{equation}
Moreover, $0\le p_0(x)-(1-D(x))\le Q(x)$. Since $D(x)\le mK/n$ for $x\in A_n^m$, and $m$, $R$, and $K$ are fixed, the expansion $\log(1-z)=-z+O(z^2)$ gives
$(n-m-R)\log(1-D(x))=-nD(x)+O(n^{-1})$.
Exponentiating, we obtain
$
(1-D(x))^{n-m-R}=e^{-nD(x)}+O(n^{-1}),
$
uniformly for $x\in A_n^m$. A telescoping argument gives~$|a^N-b^N|\le N|a-b|$ for $a,b\in[0,1]$.
Combining this with the preceding approximation yields
\[
\big|p_0(x)^{n-m-R}-e^{-nD(x)}\big|
\le nQ(x)+O(n^{-1}).
\]
Since $(n-m)_R=n^R(1+O(n^{-1}))$ and~\eqref{eq:midStep}, the preceding estimates and the triangle
inequality give, for some constant $C=C(K,\ell,m)>0$,
\[\bigg|
\frac{(n-m)_R}{n^R}p_0(x)^{n-m-R}
\prod_{i=1}^m(np_i(x))^{r_i}
-e^{-nD(x)}\prod_{i=1}^m(nd_i)^{r_i}
\bigg|
\le C\big(nQ(x)+n^{-1}\big).
\]
Dividing by $\prod_{i=1}^m r_i!$ gives the corresponding bound for the difference between the multinomial expression
in~\eqref{eqn:multinom} and the product of Poisson probabilities.
As observed above, the conditional joint degree probability
differs from~\eqref{eqn:multinom} by at most $nQ(x)$.
Thus, enlarging $C$ if necessary, we obtain
\[
\bigg|
\PP\big(d_G(i)=r_i\text{ for every }i\in[m]
       \mid X_i=x_i,\ i\in[m],\ \mathcal I\big)
-
e^{-nD(x)}\prod_{i=1}^m\frac{(nd_i)^{r_i}}{r_i!}
\bigg|
\le C\big(nQ(x)+n^{-1}\big).
\]
Integrating this error over $A_n^m$, by
\Cref{eqn:thm:mixed-factorial:Q} we have
\[ 
C\Big(
n\int_{A_n^m}Q(x)\,dx+\frac{a_n^m}{n}
\Big)
=o(n^{-m}).
\]
It remains to check that conditioning on $\mathcal I$ does not affect the
asymptotics. For $t>0$, set
\[
E_t:=\int_{S_W(t)^2}W(x,y)\,dx\,dy
\quad\text{and}\quad
h_t(z):=\int_{S_W(t)}W(x,z)\,dx.
\]
Then, $J(t)=\int_0^1h_t(z)^2\,dz$ and
$E_t=\int_{S_W(t)}h_t(y)\,dy$, and so Cauchy--Schwarz gives
$E_t^2\le g_W(t)J(t)$. Since $\lim_{t\to 0} g_W(t)/t=O(1)$ (by~\Cref{thm:measure-lowdeg}) and
$\lim_{t\to 0} J(t)/t^3=0$ (by \Cref{prop:second-order-vague}), we obtain $E_t/t^2\to 0$ as $t\to 0$. Thus, for $m\ge2$, the
integral over $A_n^m$ of the probability that at least one edge is present
between distinguished vertices is at most
\[
\binom m2a_n^{m-2}E_{K/n}=o(n^{-m}).
\]
Note that for $m=1$ there is no such edge. 
After multiplication by $(n)_m$, both the multiple-adjacency contribution
and the contribution from configurations containing an edge between two
distinguished vertices are $o(1)$. Hence, the contribution to the mixed factorial moment from
types with
$(X_1,\ldots,X_m)\in A_n^m$ is
\[
\frac{(n)_m}{n^m}
\prod_{i=1}^m
\Big[
\frac{n}{r_i!}\int_{S_W(K/n)}
(nd_W(x))^{r_i}e^{-nd_W(x)}\,dx
\Big]
+o(1).
\]
By \Cref{lem:weighted-truncated-laplace}, this converges to
\[
\prod_{r=0}^{\ell-1}
\Big(
\frac{c}{r!}\int_0^K u^re^{-u}\,du
\Big)^{m_r}.
\]
It remains to remove the restriction that all distinguished types lie in
$A_n=S_W(K/n)$. Let
$M:=\max_{i\in[m]}d_W(X_i)$.
Suppose that the distinguished vertices have respective degrees
$r_1,\ldots,r_m$, and choose one whose graphon degree equals $M$. Since the degree of the chosen vertex is at most~$\ell-1$, it has at most~$\ell-1$ neighbors
among the $n-m$ remaining vertices. Conditional on the distinguished types, let $Z$ denote the number of
neighbors of this vertex among the $n-m$ remaining vertices. Then
$Z\sim\operatorname{Bin}(n-m,M)$, and the event under consideration
requires $Z\le\ell-1$. Hence, because $m$ and $\ell$ are fixed and using that $n-m-s\ge n/2$ for all $0\le s<\ell$ and sufficiently large $n$, we have
\begin{align*}
\PP(Z\le\ell-1 \mid X_1,\ldots,X_m)
& = 
\sum_{s=0}^{\ell-1}\binom{n-m}{s}M^s(1-M)^{n-m-s} \\
& \le
\sum_{s=0}^{\ell-1}\frac{(nM)^s}{s!}e^{-nM/2} \\
& \le
e\big(1+(nM)^{\ell-1}\big)e^{-nM/2},
\end{align*}
where the last inequality follows by considering separately the cases $nM\le 1$ and $nM>1$.
Therefore the contribution from
$M>K/n$ is bounded by
\[
e(n)_m\EE\Big[
\big(1+(nM)^{\ell-1}\big)e^{-nM/2}
\Indic_{\{M>K/n\}}
\Big].
\]
Also, as the distinguished types are independent,
$\PP(M\le t)=g_W(t)^m=O(t^m)$
as $t\to 0$. The integration-by-parts argument used in the proof of
\Cref{prop:static-d} applies here with the distribution function
$g_W(t)$ replaced by $g_W(t)^m$. Indeed, writing
$h(u):=(1+u^{\ell-1})e^{-u/2}$, for all sufficiently large $K$
the function $h$ is decreasing on $[K,\infty)$.
Since $\mathbb P(M\le t)=g_W(t)^m=O(t^m)$ as $t\to0$,
Stieltjes integration by parts gives
\[\limsup_{n\to\infty}\
n^m\mathbb E\big[h(nM)\Indic_{\{M>K/n\}}\big]
\le C\int_K^\infty u^m(-h'(u))\,du,
\]
for some constant $C>0$ independent of $K$.
The integral tends to zero as $K\to\infty$ because
$h'$ decays exponentially. Since $(n)_m\le n^m$, we conclude that
\[
\lim_{K\to\infty}\limsup_{n\to\infty}\ 
(n)_m\EE\Big[
\big(1+(nM)^{\ell-1}\big)e^{-nM/2}
\Indic_{\{M>K/n\}}
\Big]
=0.
\]
We may thus first let $n\to\infty$ and then $K\to\infty$. Since
$\int_0^\infty u^re^{-u}\,du=r!$, we obtain the statement's first claim.
These are precisely the mixed factorial moments of independent
Poisson random variables with mean~$c$, so the method of moments gives the joint convergence.
\end{proof}
We are now ready to prove \Cref{thm:distribution}.
\begin{proof}[Proof of Theorem~\ref{thm:distribution}]Let $W$ be a graphon such that $\lim_{t\to 0}\frac{g_W(t)}{t}=c$ for some $c\in (0,\infty)$ and $\lim_{t\to 0}\eta(t)=0$, where $\eta(t)$ is defined as in~\eqref{def:eta}. By Corollary~\ref{cor:eta-implies-second-order}, we have that  $\kappa_n\to 0$ in the vague topology. Therefore, by \Cref{thm:measure-lowdeg}, the hypotheses of \Cref{thm:mixed-factorial} hold. Also, note that for every fixed $r\ge1$, the event~$\{\delta(G(n,W))\ge r\}$ is the same event as
$\{N^{(0)}=\cdots=N^{(r-1)}=0\}$. Therefore, by~\Cref{thm:mixed-factorial}, we have
$\PP(\delta(G(n,W))\ge r)\xrightarrow[n]{} e^{-cr}$ from which the result follows.
\end{proof}
 

\section{The $W$-random graph process}\label[section]{sec:dynamic}
Let $W$ be a graphon and $G\sim G(n,W)$. Recall that the $W$-random graph process of $G$ is the sequence of graphs $(G_k)_{k=0}^{e(G)}$, where $G_0$ is the empty graph with vertex set $[n]$ and $G_{e(G)}=G$, and, for $k\ge 1$, $G_k$ is obtained from $G_{k-1}$ by adding a random edge from $E(G)\setminus E(G_{k-1})$. Equivalently, we can consider a uniform random ordering $e_1,\ldots, e_{e(G)}$ of the edges of $G$ and define $G_k$ as the graph consisting of the first~$k$ edges of this ordering.  It will be useful to view this process as follows.
Conditional on the final graph $G$, assign independently to each edge
$e\in E(G)$ a priority $U_e\sim\operatorname{Unif}[0,1]$ and order the
edges by increasing priority. Since the priorities are almost surely distinct, ordering the edges by increasing priority produces a uniformly random ordering of $E(G)$, and hence has the same distribution as the~$W$-random graph process. 

For $0\le p\le1$, we may use the same priorities $(U_e)_{e\in E(G)}$ to define a random subgraph $H_p\subseteq G$ consisting of the edges
$e\in E(G)$ with $U_e\le p$. We call $(H_p)_{0\le p\le1}$ the
\emph{thinning process} of $G$, and note that~$H_q\subseteq H_p$ deterministically for every $q\le p$. Conditional on the vertex types $X_1,\ldots, X_n$, the edges of~$H_p$ are independent with edge probabilities $pW(X_i,X_j)$, and so
$H_p$ has the same distribution as $G(n,pW)$. Moreover, if $Y^{(p)}:=e(H_p)$, then
$H_p=G_{Y^{(p)}}$ almost surely under this coupling, as every edge has the same probability of appearing in $H_p$. Let $Y:=e(G)$ and $M_n:=\mathbb E[Y]=\binom{n}{2}\|W\|_1$, and let $p_\pm:=(1\pm\varepsilon)k/M_n$ for a fixed small $\varepsilon>0$. 
Combining \Cref{lem:Yconc} with the conditional binomial concentration of $Y^{(p)}$, we obtain that the event $Y^{(p_-)}\le k\le Y^{(p_+)}$ holds with probability tending to one, whenever $k\to\infty$ and $k=o(M_n)$. It follows that
$H_{p_-}\subseteq G_k\subseteq H_{p_+}$
deterministically on the event~$Y^{(p_-)}\le k\le Y^{(p_+)}$, and therefore we can bound the hitting times in the $W$-random graph process by studying the thinning process. We refer to this comparison as the
\emph{thinning sandwich}.

\subsection{Degree counts in the edge-revelation process}\label[section]{sec:process}
In this section, we derive the estimates that will be used in the hitting-time analysis. The main observation is that, conditional on the realization of
$G\sim G(n,W)$, the degree of any fixed vertex in $G_k$ has a
hypergeometric distribution, since the edges of $G$ are revealed
in a uniformly random order. 
Keeping this hypergeometric random variable rather than replacing it immediately by independent thinning leads to upper and lower bounds with parameters that agree to first order in the regimes of interest.
Mirroring the notation~$N^{(r)}$ from Section~\ref{sec:static}, for every integer $r\ge0$ define
\[
N_k^{(r)}:=|\{j\in[n]:d_{G_k}(j)=r\}|,
\]
where the subscript $k$ records explicitly the time in the edge-revelation process. For $\ell\ge1$, we also write
\[
N_k^{(<\ell)}:=\sum_{r=0}^{\ell-1}N_k^{(r)}
=|\{j\in[n]:d_{G_k}(j)<\ell\}|.
\]
Then, since $\tau_\ell(W)$ is the hitting time for $\delta(G_k)\geq \ell$, we have $\{\tau_\ell(W)\le k\}=\{N_k^{(<\ell)}=0\}$.

We need two auxiliary probabilistic results. The first one states that the number of edges in $G(n,W)$ concentrates around its expectation. 

\begin{lemma}[{see \cite[Corollary 10.4]{LovaszBook}}]\label[lemma]{lem:Yconc}
Let $W$ be a graphon with $\|W\|_1>0$, $Y:=e(G(n,W))$ and $M_n:=\mathbb E[Y]=\binom{n}{2}\|W\|_1$. Then, there exist a constant $c_W>0$ and an integer $n_0=n_0(W)$ such that, for all~$n\ge n_0$ and all $0<\zeta\le1$,
\[
\PP\big(|Y-M_n|>\zeta M_n\big)\le4e^{-c_W\zeta^2n}.
\]
Consequently, if $\zeta=\zeta(n)\to0$ and $\zeta^2n/\log n\to\infty$, then
$n\cdot\PP\big(|Y-M_n|>\zeta M_n\big)=o(1)$.
\end{lemma}
We also need the following lemma. 
\begin{lemma}\label[lemma]{lem:binomial}
     Let $Z\sim\mathrm{Bin}(N,q)$, let $0\le r\le N$ be an integer and let $z\in(0,1]$. Then
\[
\mathbb{E}\Big[\binom{Z}{r}z^{Z-r}\Big]=\binom{N}{r}q^r\big(1-q+qz\big)^{N-r}.
\]\end{lemma}
\begin{proof}
Indeed, write $Z=\sum_{v=1}^N\xi_v$ with $\xi_1,\dots,\xi_N$ independent $\mathrm{Bernoulli}(q)$ variables. 
Then, since on the event $\{\xi_v=1 \text{ for all } v\in S\}$ we have $Z-r=\sum_{v\notin S}\xi_v$,
\[\binom{Z}{r}z^{Z-r}=\sum_{S\subseteq[N],\,|S|=r}\ \prod_{v\in S}\xi_v\prod_{v\notin S}z^{\xi_v}.\]
By independence, each of the $\binom{N}{r}$ summands has expectation $q^r(1-q+qz)^{N-r}$.
\end{proof}

In what follows, we fix a sequence $\zeta=\zeta(n)\to0$ satisfying
$\zeta^2n/\log n\to\infty$, as in the second assertion of
\Cref{lem:Yconc}, and set $Y_\pm:=(1\pm\zeta)M_n$ and 
$\mathcal E:=\{Y_-\le Y\le Y_+\}$. The next proposition is the main result of this subsection. It bounds the quantity $\EE[N_k^{(r)}]$ from above and below by explicit expressions, similar to those in Lemma~\ref{lem:first-moment-static}.

\begin{proposition}\label[proposition]{prop:Z-two-sided}
Let $W$ be a graphon with $\|W\|_1>0$, let $G\sim G(n,W)$ and let $(G_k)_{k\ge 0}$ be the $W$-random graph process. Assume that $k=k(n)\to\infty$, $k=o(M_n)$, and $k\le Y_--n$ for all sufficiently large $n$. Fix $r\ge0$ and define\footnote{Although $a_{r,n}^+$ does not depend on $r$, we retain the subscript to keep the upper- and lower-bound notation parallel.} 
\[
a_{r,n}^+:=\frac{k}{Y_-}.
\]
Also, let
\[ 
 \theta_{r,n}^+:=\frac{k-r}{Y_+-r}\quad\text{and}\quad
 b_{r,n}^+:=1-e^{-\theta_{r,n}^+}.
\]
Set $a_{0,n}^-:=1$, define $a_{r,n}^-:=(k-r+1)/Y_+$ for $r\ge1$, and set
\[
 \sigma_{r,n}^-:=\frac{k-r}{Y_--n+1-(k-r)}
 \quad\text{and}\quad b_{r,n}^-:=1-e^{-\sigma_{r,n}^-}.
\]
Then,
\[
\EE[N_k^{(r)}]
\le n\binom{n-1}{r}(a_{r,n}^+)^r
I_r\Big((n-1-r)b_{r,n}^+\Big)
+n\cdot\PP(\mathcal E^c).
\]
Furthermore, for every $\delta=\delta(n)\in(0,1)$,
\[
\EE[N_k^{(r)}]
\ge n\binom{n-1}{r}(a_{r,n}^-)^r
\Bigg(
I_r\Big(\tfrac{(n-1-r)b_{r,n}^-}{1-b_{r,n}^-\delta}\Big)
-e^{-\frac{(n-1-r)b_{r,n}^-\delta}{1-b_{r,n}^-\delta}}
\Bigg)
-B_r n\cdot\PP(\mathcal E^c),
\]
where $B_r$ depends only on $r$.
\end{proposition}

\begin{proof}
Fix a vertex $j\in[n]$, and let $Y_j:=d_G(j)$ and
$D_j:=d_{G_k}(j)$. Conditional on the final graph $G$,
the set~$E(G_k)$ is a uniformly random $k$-subset of $E(G)$,
so $D_j$ has a hypergeometric distribution.
Then, for~$r\le Y_j$ and all sufficiently large $n$,
\[
\PP(D_j=r\mid G)
=\binom{Y_j}{r}\frac{(k)_r}{(Y)_r}
\frac{(Y-k)_{Y_j-r}}{(Y-r)_{Y_j-r}}.
\]
On the event $\mathcal E$, we have $Y_-\le Y\le Y_+$. Since $Y_j\le n-1$ and $k\le Y_--n$, the factors in the following product are nonnegative. Moreover, $Y-r-s\le Y_+-r$. Thus, we obtain
\begin{equation}\label{eqn:propZ:product}
\frac{(Y-k)_{Y_j-r}}{(Y-r)_{Y_j-r}}
=\prod_{s=0}^{Y_j-r-1}
\Big(1-\frac{k-r}{Y-r-s}\Big)
\le e^{-\theta_{r,n}^+(Y_j-r)}.
\end{equation}
On $\mathcal E$, we have
$(k)_r/(Y)_r\le(k/Y_-)^r$, and hence
$
\PP(D_j=r\mid G)
\le
\binom{Y_j}{r}(a_{r,n}^+)^r
e^{-\theta_{r,n}^+(Y_j-r)}$.
Consequently, for every $x\in[0,1]$,
\[
\PP(D_j=r,\mathcal E\mid X_j=x)
=\EE\big[
\Indic_{\mathcal E}\PP(D_j=r\mid G)
\mid X_j=x\big]
\le
(a_{r,n}^+)^r
\EE\Big[
\binom{Y_j}{r}
e^{-\theta_{r,n}^+(Y_j-r)}
\mid X_j=x\Big].
\]
Conditional on $X_j=x$, the random variable $Y_j$
has distribution $\operatorname{Bin}(n-1,d_W(x))$.
Therefore, applying \Cref{lem:binomial} with
$z=e^{-\theta_{r,n}^+}=1-b_{r,n}^+$, we obtain
\[
\PP(D_j=r,\mathcal E\mid X_j=x)
\le
\binom{n-1}{r}(a_{r,n}^+)^r d_W(x)^r
\big(1-b_{r,n}^+d_W(x)\big)^{n-1-r}.
\]
Using $1-u\le e^{-u}$, integrating over $x\in[0,1]$,
and recalling the definition of $I_r$ in
\Cref{sec:static}, we conclude that
\[
\PP(D_j=r,\mathcal E)
\le
\binom{n-1}{r}(a_{r,n}^+)^r
I_r\Big((n-1-r)b_{r,n}^+\Big).
\]
Summing over $j$ and adding the contribution of $\mathcal E^c$, which is at most
$n\cdot\PP(\mathcal E^c)$, proves the upper bound for~$\EE[N_k^{(r)}]$.

For the lower bound, on $\mathcal E$ we have $(k)_r/(Y)_r\ge(a_{r,n}^-)^r$. For each factor in~\eqref{eqn:propZ:product}, put $u_s:=(k-r)/(Y-r-s)$. Since $Y_j\le n-1$ and $k\le Y_--n$, one has $u_s<1$. Using $1-u\ge e^{-u/(1-u)}$ gives
\[
\frac{(Y-k)_{Y_j-r}}{(Y-r)_{Y_j-r}}
\ge e^{-\sigma_{r,n}^-(Y_j-r)}.
\]
It follows that, on $\mathcal E$,
\[
\PP(D_j=r\mid G)\ge F_j(G)
\quad\text{where}\quad
F_j(G):=\binom{Y_j}{r}(a_{r,n}^-)^r
e^{-\sigma_{r,n}^-(Y_j-r)}.
\]
Hence, $\PP(D_j=r)\ge \EE[\Indic_{\mathcal E}F_j(G)]$.
When $k/M_n=o(1)$ and $k\to\infty$, the random variable
$F_j(G)$ is bounded by a constant depending only on $r$.
Indeed, for $r=0$, we have
$F_j(G)=e^{-\sigma_{0,n}^-Y_j}\le1$.
For $r\ge1$, using~$\binom yr\le y^r/r!$ and the fact that
the maximum of $y^re^{-\sigma y}$ is $(r/(e\sigma))^r$,
the boundedness follows because $\sigma_{r,n}^-\to0$ and~$a_{r,n}^-/\sigma_{r,n}^-\to 1$.
Hence, for some constant $B_r>0$ depending only on $r$,
\[
\PP(D_j=r)
\ge \EE[F_j(G)]-B_r\PP(\mathcal E^c).
\]
Now, conditional on $X_j=x$, the random variable $Y_j$
has distribution $\mathrm{Bin}(n-1,d_W(x))$.
Applying \Cref{lem:binomial} with
$z=e^{-\sigma_{r,n}^-}=1-b_{r,n}^-$ gives
\[
\EE[F_j(G)\mid X_j=x]
=
\binom{n-1}{r}(a_{r,n}^-)^r d_W(x)^r
\big(1-b_{r,n}^-d_W(x)\big)^{n-1-r}.
\]
Integrating over $x$ and summing over $j$, it remains to
bound the resulting integral from below.

Finally, on the set $S_W(\delta)$, the inequality
$\log(1-bd_W(x))\ge-bd_W(x)/(1-b\delta)$
gives a lower bound in terms of the exponential integral
over $S_W(\delta)$. Extending that integral to $[0,1]$
adds at most
$\exp(-(n-1-r)b\delta/(1-b\delta))$. Hence, subtracting this contribution proves the stated lower bound.
\end{proof}

The next statement shows that
all of the hypergeometric parameters appearing in \Cref{prop:Z-two-sided} have the same
first-order asymptotic values.

\begin{corollary}\label[corollary]{cor:Z-parameters}Under the same hypotheses as in \Cref{prop:Z-two-sided} the following hold. For every fixed $r\ge0$,
\[
 a_{r,n}^+=\frac{k}{M_n}\Big(1+o(1)\Big),\qquad
 a_{r,n}^-=\frac{k}{M_n}\Big(1+o(1)\Big),
\]
and
\[
 (n-1-r)b_{r,n}^+=\frac{nk}{M_n}\Big(1+o(1)\Big),\qquad
 (n-1-r)b_{r,n}^-=\frac{nk}{M_n}\Big(1+o(1)\Big).
\]
Moreover, $n\cdot\PP(\mathcal E^c)=o(1)$.
\end{corollary}

\begin{proof}
The first two estimates use $r=O(1)$ and $\zeta=o(1)$. Since $k/M_n=o(1)$, both $\theta_{r,n}^+$ and $\sigma_{r,n}^-$ are~$k/M_n\cdot(1+o(1))$. The expansion $1-e^{-u}=u\big(1+O(u)\big)$ then gives the two estimates for $b_{r,n}^\pm$. The last assertion follows from \Cref{lem:Yconc} as $\|W\|_1>0$.
\end{proof}

The asymptotic formulas for $a_{r,n}^{\pm}$ and
$(n-1-r)b_{r,n}^{\pm}$ in \Cref{cor:Z-parameters} allow us to evaluate the upper and lower
bounds in \Cref{prop:Z-two-sided} explicitly in the intermediate regime where the lower
tail of the graphon degree distribution is asymptotically linear.

\subsection{A sparse analogue of Theorem~\ref{thm:hladky}}\label[section]{sec:criteria}

In this section, we extend \Cref{thm:hladky} for a rescaled graphon of the form $p_nW$, where $p_n\in [0,1]$ is a sequence such that $p_nn\to\infty$. For a graphon $W$ and a sequence $(a_n)_{n\in\mathbb N}$ with $a_n>0$ for all $n\in\mathbb N$, let~$R_W(a_n):=\sup_{0<t\le a_n}g_W(t)/t$.
\begin{theorem}\label[theorem]{thm:triangular}Let $W$ be a graphon, let $p_n\in(0,1]$ satisfy $np_n\to\infty$, and let $H_n\sim G(n,p_nW)$.
\begin{enumerate}[label=(\alph*)]
\item\label{thm:triangular:itm1} If there exists $q_n\to0$ such that
$n g_W(q_n/(np_n))\to\infty$, then a.a.s.\ $H_n$ contains an isolated vertex.

\item\label{thm:triangular:itm2} If there exists $a_n\to0$ such that
$np_na_n/\log n\to\infty$ and $R_W(a_n)/p_n\to0$, then $\delta(H_n)\to\infty$ in probability.

\item\label{thm:triangular:itm3} If, for some fixed $q>0$,
$\liminf_{n\to\infty}n g_W(q/(np_n))>0$, then
  \[
  \liminf_{n\to\infty}\PP(\delta(H_n)=0)>0.
  \]
\item\label{thm:triangular:itm4} If there are $a_n\to0$ and $\Psi>0$ such that
$np_na_n/\log n\to\infty$ and $R_W(a_n)\le\Psi p_n$ for all sufficiently large $n$, then, for every fixed $\ell$,
\[\liminf_{n\to\infty}\PP(\delta(H_n)>\ell)>0.
\]
\end{enumerate}
\end{theorem}

\begin{proof}
The proofs rely heavily on the first- and second-moment arguments developed in \Cref{sec:static}. While parts~\ref{thm:triangular:itm1} and~\ref{thm:triangular:itm3} adapt second-moment arguments, parts (b) and (d) extend the first-moment arguments for the rescaled graphon $p_nW$.

Observe that $d_{p_nW}(x)=p_nd_W(x)$ for every $x\in [0,1]$, so $g_{p_nW}(t)=g_W(t/p_n)$.

For part~\ref{thm:triangular:itm1}, choose $c_n\to\infty$ with $c_n\le n g_W(q_n/(np_n))$ and let $S_n\subseteq S_W(q_n/(np_n))$ have measure $c_n/n$. Let $V_n$ be the number of isolated vertices whose types lie in $S_n$. Then, in the graphon $p_nW$, every such type has graphon degree at most $q_n/n$. The second-moment bound recorded in the proof of \Cref{prop:static-ac}, with $Lq$ replaced by $c_n$ and $q$ replaced by $q_n$, 
implies that $\PP(V_n>0)\to1$, since $\EE[V_n]\ge c_n(1-q_n/n)^{n-1}=(1+o(1))c_n$ and $\EE[V_n^2]\le\EE[V_n]+c_n^2$. 
Part~\ref{thm:triangular:itm3} follows from the same second-moment bound: choose a fixed $c>0$ satisfying
$c<\liminf_{n\to\infty}
n g_W\big(q/(np_n)\big)$,
and, for all sufficiently large $n$, choose
$S_n\subseteq S_W(q/(np_n))$ with $\lambda(S_n)=c/n$.

For part~\ref{thm:triangular:itm2}, apply the exact first-moment identity of \Cref{lem:first-moment-static} to the graphon $p_nW$. Since $d_{p_nW}(x)=p_nd_W(x)$, the upper bound for the number of vertices of degree $r$ has the form
$n\binom{n-1}{r}p_n^r I_r(C)$, with $C:=(n-1-r)p_n$. Repeating the proof of \Cref{lem:Ir-small-tail}, but using the uniform bound $g_W(t)\le R_W(a_n)t$ on~$t\le a_n$, gives, for every fixed $r$,
$I_r(C)\le \Lambda_r C^{-r}(e^{-Ca_n/2}+R_W(a_n)/C)$. Multiplying each of the two terms on the right-hand side of
the preceding inequality by $n\binom{n-1}{r}p_n^r$ gives contributions bounded by $O(ne^{-Ca_n/2})=O(ne^{-np_na_n/2})$ and $O(R_W(a_n)/p_n)$, respectively, both of which tend to zero by the hypotheses of part~\ref{thm:triangular:itm2}. Hence, for every fixed $\ell$, the expected number of vertices of degree below $\ell$ tends to zero. By Markov's inequality, a.a.s.\ there are no such vertices, which proves that $\delta(H_n)\to\infty$ in probability.

For part~\ref{thm:triangular:itm4}, apply the proof of
\Cref{prop:static-d} to the graphon $p_nW$.
Its lower-tail function is~$g_{p_nW}(u)=g_W(u/p_n)$,
so the assumption $R_W(a_n)\le\Psi p_n$ gives
$g_{p_nW}(u)\le\Psi u$ for $u\le p_na_n$. Fix a constant $c>0$, to be chosen sufficiently large
depending only on $\Psi$ and $\ell$, and let $S_n$ be
the exceptional set given by \Cref{lem:exceptional}
for the graphon $p_nW$. Thus, $S_n$ contains all types
of degree at most $2c/n$, and every type outside $S_n$
retains at least one quarter of its graphon degree
outside $S_n$.
Let $\mathcal S_n$ be the event that none of the
sampled vertex types belong to $S_n$.
Conditional on $\mathcal S_n$, the expected number of vertices of degree $r$ is bounded by a constant
multiple of
\[n^{r+1}\int_{\{d_{p_nW}>2c/n\}}
d_{p_nW}(x)^r e^{-n d_{p_nW}(x)/8}\,dx.
\]
Splitting this integration domain as $\{2c/n<d_{p_nW}(x)\le p_na_n\}\cup \{d_{p_nW}(x)>p_na_n\}$,
the contribution from~$d_{p_nW}(x)>p_na_n$,
including the factor $n^{r+1}$, is at most
$n^{r+1}e^{-np_na_n/8}=o(1)$ as
$np_na_n/\log n\to\infty$.
On the range $2c/n<d_{p_nW}(x)\le p_na_n$,
the estimate from the proof
of \Cref{prop:static-d} (with $t_0=p_na_n$) applies uniformly in $n$,
with constants depending only on $\Psi$ and $\ell$. Therefore, we may choose a truncation constant $c$ sufficiently large to get the desired lower bound.
\end{proof}

\subsection{Connectivity and the last isolated vertex}\label[section]{sec:connectivity-hitting}
In this subsection, we prove \Cref{theorem:hitting-time}, which is the analogue for the $W$-random graph process of the classical hitting-time theorem asserting that the random graph process becomes connected exactly when its last isolated vertex disappears (see~\cite{BollobasThomason1985}). The proof of \Cref{theorem:hitting-time} builds on the strategy of the classical homogeneous case. The main idea is to construct two thinning parameters, say $p_n^- < p_n^+$, so that a.a.s.\ $H_{p_n^-}$ consists of a large connected component and a few isolated vertices (with at least one remaining), whereas~$H_{p_n^+}$ is a.a.s.\ connected. The thinning sandwich therefore gives lower
and upper bounds on the hitting times for connectivity and minimum degree one. To show that both hitting times coincide, we show that when $(1+\varepsilon)e(H_{p_n^-})\le k<\tau_{\mathrm{conn}}(W)$, a.a.s.\ an isolated vertex in $G_k$ can only form edges with vertices in the largest component. Therefore, the moment in which the last isolated vertex disappears is the same as when the graph becomes connected. Though this strategy mirrors the proof in the classical setting, we need to introduce new technical tools for the inhomogeneous case
(for example, we use the Laplace transform to establish sharp upper bounds for the expected number of isolated vertices).

The first result we need is a uniform factorial-moment estimate for the number of isolated vertices in~$H_p$, expressed in terms of the Laplace transform $\mathcal L(C)=I_0(C)=\int_0^1 e^{-Cd_W(x)}\,dx$ introduced in \Cref{ssec:alternativeProof}. Unlike \Cref{thm:mixed-factorial}, this estimate does not require the second-order vague decorrelation condition. Indeed, thinning by $p$ gives $d_{pW}(x)=p\,d_W(x)$ and $q_{pW}(x,y)=p^2q_W(x,y)$, so
\[
q_{pW}(x,y)\le\frac p2\big(d_{pW}(x)+d_{pW}(y)\big).
\]
Thus, as $p\to0$, pairwise neighborhood overlaps become uniformly negligible relative to the corresponding degrees. We exploit this bound directly to control joint isolation probabilities, without invoking vague convergence of the second-order measures. The following lemma makes this precise.

\begin{lemma}[Isolated vertices under thinning]\label[lemma]{lem:thinned-isolated-factorial}
Let $W$ be a connected graphon with $\lim_{t\to 0}\frac{g_W(t)}t=0$. Let~$C_n\to\infty$ and $p_n:=C_n/n$, and suppose that $p_n\to0$. Let $N^{(0)}_{p_n}$ denote the number of isolated vertices in~$H_{p_n}$. If $\mu_n:=n\mathcal L(C_n)$ satisfies
$1\le\mu_n\le\log n$, then, for every integer sequence $J_n=O(\log(1/p_n))$,
\[
\EE[(N^{(0)}_{p_n})_j]
=\big(1+o(1)\big)\mu_n^j
\]
uniformly for $1\le j\le J_n$.
\end{lemma}

\begin{proof}
Fix $1\le j\le J_n$. By exchangeability, label vertices
$1,\ldots,j$ as the distinguished vertices and vertices
$j+1,\ldots,n$ as the remaining vertices, and condition on the
distinguished types $X_i=x_i$, $i\in[j]$. 
For any of the remaining vertices, the probability of avoiding all distinguished
vertices in $H_{p_n}$ is
\[
R:=\int_0^1
\prod_{i=1}^j\big(1-p_nW(x_i,u)\big)\,du.
\]
Write $d_i:=d_W(x_i)$ and $D:=\sum_{i=1}^j d_i$. The union bound and the second Bonferroni inequality give
\[
1-p_nD
\le R\le
1-p_nD+
p_n^2\sum_{1\le i<r\le j}q_W(x_i,x_r).
\]
Since
$q_W(x_i,x_r)
\le \min\{d_i,d_r\}
\le (d_i+d_r)/2$,
we have
$\sum_{1\le i<r\le j}q_W(x_i,x_r)
\le (j-1)D/2$.
Consequently,
\[
1-p_nD
\le R\le
1-p_n\Big(1-\frac{p_n(j-1)}2\Big)D.
\]
Conditional on the types of the distinguished vertices, 
\[\mathbb P(1,\ldots,j\text{ are isolated})=
R^{n-j}
\prod_{1\le i<r\le j}
\big(1-p_nW(x_i,x_r)\big)
\]
and, since $0\le W\le 1$, this value lies between $R^{n-j}(1-p_n)^{\binom{j}{2}}$ and $R^{n-j}$. Moreover, $p_nj=o(1)$ uniformly for~$j\le J_n$ (because, since $1/p_n\to\infty$ by hypothesis, $p_nJ_n=O(\log(1/p_n)/(1/p_n))=o(1)$). Hence, using that~$\log(1-u)\ge-u/(1-u)$ for $u\in (0,1)$, we have
\[
R^{n-j}
\ge
(1-p_nD)^{n-j}
\ge
\exp\Big(
-\frac{(n-j)p_nD}{1-p_nD}
\Big)
\ge e^{-A_{j,n}D},
\]
where
$A_{j,n}:=\frac{(n-j)p_n}{1-p_nj}$ (using $D\le j$).
For the upper bound,
\[
R^{n-j}
\le
\exp\Big(
-(n-j)p_n
\big(1-\frac{p_n(j-1)}2\big)D
\Big)
=e^{-B_{j,n}D},
\]
where
$B_{j,n}:=
(n-j)p_n\big(1-\frac{p_n(j-1)}2\big)$.
Integrating over $x_1,\ldots,x_j$ and using
$D=\sum_i d_i$ therefore gives
\begin{equation}\label{eq:joint-isolation-squeeze}
(1-p_n)^{\binom j2}\mathcal L(A_{j,n})^j
\le
\PP(1,\ldots,j\text{ are isolated})
\le
\mathcal L(B_{j,n})^j.
\end{equation}
We next show that replacing $C_n$ by either
$A_{j,n}$ or $B_{j,n}$ changes the corresponding
$j$-th powers of $\mathcal L$ by a multiplicative factor
$1+o(1)$, uniformly for $1\le j\le J_n$. Note first that, since $C_n=np_n$,
$A_{j,n}/C_n
=1+O\big(p_nj+\frac{j}{n}\big)$ and
$B_{j,n}/C_n
=1+O\big(p_nj+\frac{j}{n}\big)$
uniformly for $1\le j\le J_n$ (where we use again that $p_nJ_n=o(1)$). Put
$\varepsilon_n:=
p_nJ_n+\frac{J_n}{n}=o(1)$. To prove that
$\mathcal L(A_{j,n})^j=(1+o(1))\mathcal L(C_n)^j$ and
$\mathcal L(B_{j,n})^j=(1+o(1))\mathcal L(C_n)^j$, it suffices to show that, uniformly over all $\gamma$ satisfying
$|\gamma|\le C\varepsilon_n$ for any fixed $C<\infty$, we have
\begin{equation}\label{eq:L-uniform-stability}
\frac{\mathcal L((1+\gamma)C_n)}{\mathcal L(C_n)}
=
1+o\Big(\frac1{J_n}\Big).
\end{equation}
Suppose first that
$\operatorname*{ess\,inf}d_W>0$. Choose $\delta_0>0$ such that
$d_W(x)\ge\delta_0$ almost everywhere. Since
$\mathcal L(C_n)\le e^{-\delta_0C_n}$
and $n\mathcal L(C_n)=\mu_n\ge1$, we have $C_n=O(\log n)$. Also,
$J_n=O(\log(1/p_n))=O(\log n)$. Hence
\[
J_n\varepsilon_nC_n
=
O\Big(
\frac{C_n^2J_n^2}{n}
+
\frac{C_nJ_n^2}{n}
\Big)
=o(1).
\]
Note that for $|\gamma|\le C\varepsilon_n$ and $0\le d_W\le1$, we have
\[
e^{-|\gamma|C_n}
\le
\frac{e^{-(1+\gamma)C_nd_W(x)}}
     {e^{-C_nd_W(x)}}
\le
e^{|\gamma|C_n}.
\]
Thus, after integration,
\[
e^{-|\gamma|C_n}
\le
\frac{\mathcal L((1+\gamma)C_n)}{\mathcal L(C_n)}
\le
e^{|\gamma|C_n}.
\]
Since $J_n|\gamma|C_n=o(1)$ uniformly, this proves
\eqref{eq:L-uniform-stability} when $\operatorname*{ess\,inf} d_W>0$.

Suppose now that
$\operatorname*{ess\,inf}d_W=0$. The bound $\mu_n\le\log n$ implies
$\frac{C_n}{\log n}\longrightarrow\infty$.
Indeed, otherwise, along a subsequence,
$C_n\le A\log n$ for some constant $A>0$. Since the essential infimum is
zero, we may choose a fixed $a>0$ such that $Aa<1/2$ and
$g_W(a)>0$. Then,
$
n\mathcal L(C_n)
\ge
ng_W(a)e^{-C_na}
\ge
g_W(a)n^{1/2}$,
contradicting that $\mu_n\le\log n$.

Set
$\Phi_n:=2\log(1/p_n)$, and note that, since $\Phi_n\le2\log n$, we have $\Phi_n/C_n\to0$. Choose $a>0$ small enough so 
that $g_W(t)\le t$ for $0<t\le a$. Then, integration by parts gives
\begin{align*}
T_n:=
\int_{\{C_nd_W>\Phi_n\}}
e^{-C_nd_W(x)}\,dx 
& = \int_{\{\frac{\Phi_n}{C_n}< d_W\le a\}}e^{-C_nd_W(x)}\,dx
+ \int_{\{d_W> a\}}e^{-C_nd_W(x)}\,dx
\\ &
\le e^{-C_na}
+C_n\int_{\Phi_n/C_n}^a g_W(t)e^{-C_nt}\,dt
\\
& \le
2e^{-C_na}
+
\frac1{C_n}\int_{\Phi_n}^{\infty}ue^{-u}\,du,
\end{align*}
and therefore
$T_n
\le
2e^{-C_na}
+
(\Phi_n+1)e^{-\Phi_n}/C_n$.
Since $e^{-\Phi_n}=p_n^2$, 
$\mathcal L(C_n)=\frac{\mu_n}{n}
=\frac{\mu_np_n}{C_n}$ and $\mu_n\ge 1$,
we obtain
\[
\frac{J_nT_n}{\mathcal L(C_n)}
\le2J_nne^{-C_na}
+
J_n(\Phi_n+1)p_n
=o(1).
\]
Using that $J_n=O(\Phi_n)$ and $\Phi_n\le 2\log n$,  the first term tends to zero because
$C_n/\log n\to\infty$, while the second tends to zero because $1/p_n\to\infty$ and $p_n\Phi_n^2=O((\log(1/p_n)/\sqrt{1/p_n})^2)=o(1)$. 

Now let $|\gamma|\le C\varepsilon_n$. Since $J_n=O(\Phi_n)$, we have
\[
J_n|\gamma|\Phi_n
=
O\Big(
p_nJ_n^2\Phi_n+\frac{J_n^2\Phi_n}{n}
\Big)
=o(1),
\]
as the two terms are respectively
$O(p_n\Phi_n^3)$ and $O(\Phi_n^3/n)$. 

On the set $\{C_nd_W\le \Phi_n\}$,
$e^{-(1+\gamma)C_nd_W}
=
e^{-C_nd_W}e^{-\gamma C_nd_W}$,
and hence
$e^{-|\gamma|\Phi_n}
\le
e^{-\gamma C_nd_W}
\le
e^{|\gamma|\Phi_n}$. Thus, uniformly for $|\gamma|\le C\varepsilon_n$, the contribution $I'_n$ of this set to $\mathcal L((1+\gamma)C_n)$ differs from its contribution $I_n$ to~$\mathcal L(C_n)$ by a multiplicative  factor
$1+o(1/J_n)$.

It remains to bound the contribution $T'_n$ of the set $\{C_nd_W>\Phi_n\}$ to $\mathcal L((1+\gamma)C_n)$. Put
$C_n':=(1+\gamma)C_n$ and $\Phi_n':=(1+\gamma)\Phi_n$, and note that $\{C_nd_W>\Phi_n\}=\{C_n'd_W>\Phi_n'\}$.
Applying the preceding tail estimate with $C_n'$ and $\Phi_n'$ gives
\[T'_n:=
\int_{\{C_nd_W>\Phi_n\}}
e^{-C_n'd_W(x)}\,dx
\le
2e^{-C_n'a}
+
\frac{(\Phi_n'+1)e^{-\Phi_n'}}{C_n'}.
\]
Since $\gamma=o(1)$ and
$J_n|\gamma|\Phi_n=o(1)$, we have
\[
C_n'=(1+o(1))C_n,
\quad\text{and}\quad
e^{-\Phi_n'}
=
e^{-\Phi_n}e^{-\gamma \Phi_n}
=
p_n^2\Big(1+o\Big(\frac1{J_n}\Big)\Big).
\]
Consequently, uniformly for $|\gamma|\le C\varepsilon_n$, we get $T'_n=o(\mathcal{L}(C_n)/J_n)$. Together with the corresponding estimate for $I_n$, $I'_n$ and $T_n$, this gives
\[
|\mathcal L(C_n')-\mathcal L(C_n)|
\le |I_n'-I_n|+T_n'+T_n\\
=o\Big(\frac{\mathcal L(C_n)}{J_n}\Big).
\]
Dividing by $\mathcal L(C_n)$
proves \eqref{eq:L-uniform-stability} also when
$\operatorname*{ess\,inf}d_W=0$. We may now apply \eqref{eq:L-uniform-stability} to
$A_{j,n}$ and $B_{j,n}$. Indeed, as
\[
\frac{\mathcal L(A_{j,n})}{\mathcal L(C_n)}
=
1+o\Big(\frac1{J_n}\Big)
\quad\text{and}\quad
\frac{\mathcal L(B_{j,n})}{\mathcal L(C_n)}
=
1+o\Big(\frac1{J_n}\Big)
\]
uniformly for $j\le J_n$, it follows that
$\mathcal L(A_{j,n})^j
=(1+o(1))
\mathcal L(C_n)^j$
and
$\mathcal L(B_{j,n})^j
=(1+o(1))
\mathcal L(C_n)^j$
uniformly for $1\le j\le J_n$. Finally, as
$p_nJ_n^2
=
O\big(p_n\log^2(1/p_n)\big)
=o(1)$
and 
$\frac{J_n^2}{n}=o(1)$, we have, uniformly for $1\le j\le J_n$, that
$(1-p_n)^{\binom j2}=1+o(1)$ and
$\frac{(n)_j}{n^j}=1+o(1)$. Multiplying \eqref{eq:joint-isolation-squeeze} by $(n)_j$ therefore
gives
\[
\EE[(N^{(0)}_{p_n})_j]
=
(1+o(1))
\big(n\mathcal L(C_n)\big)^j
=
(1+o(1))\mu_n^j
\]
uniformly for $1\le j\le J_n$, as required.
\end{proof}
We also need the following connectivity theorem of Hladk\'y and Viswanathan, which we will reprove in \Cref{app:static-simplifications} using the tools of Section~\ref{sec:analytic}. 

\begin{theorem}[{\cite[Theorem 4]{HV}}]\label[theorem]{thm:static-connectivity}
Suppose that \(W\) is a connected graphon such that
$\lim_{t\to 0}\frac{g_W(t)}t=0$.
Then, for every fixed \(K\in\mathbb N\),  \(G(n,W)\) is a.a.s.
\(K\)-connected.
\end{theorem}

The following lemma parallels a key step in the classical proof that, in the Erd\H{o}s--R\'enyi random graph process, connectivity occurs exactly when the last isolated vertex disappears, as established by Bollob\'as and Thomason~\cite{BollobasThomason1985}. It identifies two thinning parameters such that, at the first, the graph consists a.a.s.\ of a single non-trivial connected component and a growing but small number of isolated vertices, while at the second, the graph is connected a.a.s. As in the classical argument, it will then suffice to show that no edge joining two of these isolated vertices appears between the two parameters. Unlike in the homogeneous setting, the parameters are chosen using the degree distribution of $W$, rather than an explicit threshold. The proof uses \Cref{thm:static-connectivity} with $K=1$.
\begin{lemma}[A pre-connectivity window]
\label[lemma]{lem:preconnectivity-window}
Suppose that $W$ is a connected graphon such that
$\lim_{t\to0}\frac{g_W(t)}t=0$.
Then, there exist sequences $p_n^-<p_n^+$ and $h_n$ such that
\[p_n^-\xrightarrow[n\to\infty]{}0,
\qquad p_n^+\xrightarrow[n\to\infty]{} 0,
\qquad h_n\xrightarrow[n\to\infty]{}\infty,
\qquad h_n=o(\log n),
\qquad p_n^+h_n^2\xrightarrow[n\to\infty]{}0,
\]
and the following hold. If $\mathcal C_j^-$ denotes the number of components in $H_{p_n^-}$ with exactly $j$ vertices, then
\begin{enumerate}[label=(\roman*)]
\item\label{lem:preconnectivity-window:itm1} $\mathcal C_1^-/h_n\to1$ in probability,
\item\label{lem:preconnectivity-window:itm2} $\displaystyle \EE\Big[\textstyle\sum_{j=2}^{\lfloor n/2\rfloor}\mathcal C_j^-\Big]=o(1)$, and
\item\label{lem:preconnectivity-window:itm3} $H_{p_n^+}$ is connected a.a.s.
\end{enumerate}
In particular,
$\PP(\mathcal C_1^->2\log n)=o(1)$.
\end{lemma}

\begin{proof}
We first construct the sequence $(q_n)_{n\in\mathbb N}$ such that $H_{q_n}$ is connected a.a.s. Note that for every fixed~$p>0$, the graphon $pW$ is connected because $W$ is connected, and, moreover, $d_{pW}(x)=p\,d_W(x)$ and hence
$g_{pW}(t)=g_W(t/p)$. Therefore,
\[\frac{g_{pW}(t)}{t}
=
\frac1p\frac{g_W(t/p)}{t/p}
\xrightarrow[t\to 0]{}0.
\]
Thus, \Cref{thm:static-connectivity}, applied to $pW$ and $K=1$, shows that $G(n,pW)$ is connected a.a.s. By a standard diagonal argument we may find a sequence $(q_n)_{n\in\mathbb N}$ such that $q_n\to0$ and $H_{q_n}$ is a.a.s.\ connected. 

Recall that $W$ being connected implies that $d_W>0$ almost everywhere.  Hence, $\mathcal{L}(\cdot)$ is continuous and strictly decreasing from $1$ to $0$. For an integer $m\ge1$ and $n>m$, let $C_{n,m}$ be the unique value satisfying~$n\mathcal L(C_{n,m})=m$, and put $p_{n,m}:=C_{n,m}/n$. For fixed $m$, note that $C_{n,m}\to\infty$ as $n$ goes to infinity. Then, by \Cref{lem:Ir-small-tail}, we have
\[
p_{n,m}m^2
=mC_{n,m}\mathcal L(C_{n,m})\xrightarrow[n\to \infty]{} 0.
\]
Choose a sequence of positive integers $(h_n)_{n\in\mathbb N}$ such that $h_n\to\infty$ sufficiently slowly and
\[
h_n=o(\log n),
\qquad q_nh_n^2\to0,
\qquad p_{n,h_n}h_n^2\to0.
\]
Such a choice follows by a standard diagonal construction: keep $h_n=m$ until all three displayed requirements hold with error at most $1/m$, and then increase $m$ by one.
Set
\[
p_n^-:=p_{n,h_n},
\qquad
C_n:=np_n^-,
\quad\text{and}\quad
p_n^+:=\max\{2p_n^-,q_n\},
\]
and note that $p_n^\pm\to0$ and $p_n^+h_n^2\to0$. Since $p_n^+\ge q_n$, monotonicity implies that $H_{p_n^+}$ is connected a.a.s., so~\ref{lem:preconnectivity-window:itm3} holds. Note that by construction,
$n\mathcal L(C_n)=n\mathcal L(C_{n,h_n})=h_n$, and applying \Cref{lem:thinned-isolated-factorial} with~$j=1,2$ gives
$\EE[\mathcal C_1^-]=(1+o(1))h_n$ and
$\EE[(\mathcal C_1^-)_2]=(1+o(1))h_n^2$. Then 
$\operatorname{Var}(\mathcal C_1^-)=o(h_n^2)$, and hence (i) follows from a second-moment argument.

It remains to prove~\ref{lem:preconnectivity-window:itm2}. First, letting
$a_n:=e\,C_n\mathcal L(C_n/3)$, it follows from \Cref{lem:Ir-small-tail} that $a_n\to0$. Then, defining
\[
J_n:=\max\Big\{2,
\Big\lceil\frac{2\log(1/p_n^-)}{\log(1/a_n)}\Big\rceil\Big\},
\]
for all sufficiently large $n$, we obtain that $J_n=O(\log(1/p_n^-))$ and $a_n^{J_n}\le(p_n^-)^2$.

We now divide the possible sizes of non-trivial components of $H_{p_n^-}$ into four ranges. For $2\le j\le J_n$, compare the event that a fixed $j$-set $A$ is a component with the event that all vertices of $A$ are isolated. Conditional on the types, the requirement that no edge joins $A$ to its complement is the same in both events. By Cayley's formula, the probability that $H_{p_n^-}[A]$ is connected is at most $j^{j-2}(p_n^-)^{j-1}$, while the probability that it has no internal edge is at least $(1-p_n^-)^{\binom j2}$. Consequently,
\[
\EE[\mathcal C_j^-]
\le
\frac{j^{j-2}}{j!}(p_n^-)^{j-1}
(1-p_n^-)^{-\binom j2}
\EE[(\mathcal C_1^-)_j].
\]
Therefore, by \Cref{lem:thinned-isolated-factorial}, we have that 
\[
\EE[\mathcal C_j^-]
\le
\big(1+o(1)\big)
\frac{j^{j-2}}{j!}(p_n^-)^{j-1}h_n^j,
\]
uniformly for $2\le j\le J_n$. Since $j^{j-2}/j!\le e^j$ and $p_n^-h_n\to0$, summing over $2\le j\le J_n$ gives
\begin{equation}\label{eq:small-components-small-range}
\sum_{j=2}^{J_n}\EE[\mathcal C_j^-]
=O(p_n^-h_n^2)=o(1).
\end{equation}
We next consider the range $J_n<j\le1/p_n^-$. For a fixed $j$-set $A$, condition on the types $(x_i)_{i\in A}$ and put~$S_A:=\sum_{i\in A}d_W(x_i)$. For a given vertex outside $A$ with type $X$, the inequality
$1-\prod_{i\in A}(1-p_n^-W(x_i,X))\ge s/(1+s)$ with
$s=p_n^-\sum_{i\in A}W(x_i,X)\le p_n^-j$, followed by averaging over
its type $X$, shows that the probability of having no neighbor in $A$ is at
most $1-p_n^-S_A/(1+p_n^-j)$. Hence, the probability that no edge joins $A$ to its complement is at most
$\exp(-\sigma_{n,j}S_A)$, where
\[
\sigma_{n,j}:=
\frac{(n-j)p_n^-}{1+p_n^-j}.
\]
After averaging over the types, and again using Cayley's formula for the number of spanning trees on $A$, we have
\[
\EE[\mathcal C_j^-]
\le
\binom nj j^{j-2}(p_n^-)^{j-1}\mathcal L(\sigma_{n,j})^j.
\]
Since $j\le1/p_n^-$ and $C_n\to\infty$, one has $\sigma_{n,j}\ge C_n/3$ for all sufficiently large $n$. Therefore,
\[
\EE[\mathcal C_j^-]
\le
\frac{1}{p_n^-j^2}
\big(e\,C_n\mathcal L(C_n/3)\big)^j
=
\frac{a_n^j}{p_n^-j^2}.
\]
Recalling that $a_n^{J_n}\le(p_n^-)^2$ and $a_n\to0$, we get
\begin{equation}\label{eq:small-components-middle-range}
\sum_{J_n<j\le1/p_n^-}\EE[\mathcal C_j^-]
\le\frac1{p_n^-}\sum_{j>J_n}a_n^j
=\frac{a_n^{J_n+1}}{p_n^-(1-a_n)}
\le\frac{p_n^-a_n}{1-a_n}
=o(1).
\end{equation}
For the next range, choose $\varepsilon_0>0$ sufficiently small so that $4e\varepsilon_0<1$. By \Cref{lem:Ir-small-tail}, there is $C_0$ such that~$\mathcal L(C)\le\varepsilon_0/C$ for all $C\ge C_0$. Fix $\gamma\in(0,1/2)$ sufficiently small that $1/(4\gamma)\ge C_0$. For $1/p_n^-\le j\le\gamma n$, we have $n-j\ge n/2$ and $1+p_n^-j\le2p_n^-j$.
Hence, $\sigma_{n,j}\ge n/(4j)\ge C_0$. Then we have
\[
\EE[\mathcal C_j^-]
\le
\binom nj\mathcal L(\sigma_{n,j})^j\le \Big(\frac{en}{j}\Big)^j(\varepsilon_0/\sigma_{n,j})^j
\le
(4e\varepsilon_0)^j,
\]
and since $1/p_n^-\to\infty$,
\begin{equation}\label{eq:small-components-large-micro}
\sum_{1/p_n^-\le j\le\gamma n}\EE[\mathcal C_j^-]=o(1).
\end{equation}
Finally, consider the range $\gamma n\le j\le n/2$. By \cite[Proposition~11]{HV}, which follows from Janson's characterization of connected graph limits~\cite[Theorem~1.10(i)]{JansonConnected}, there exists $\rho=\rho(W,\gamma)>0$ such that every graph $H$ sufficiently close to $W$ in cut distance has at least $\rho |V(H)|^2$ edges across every cut whose smaller side contains at least $\gamma |V(H)|$ vertices. Since $H_1=G(n,W)$ converges to $W$ in cut distance a.a.s.~\cite[Lemma~10.16]{LovaszBook}, a.a.s.\ every cut $(A,[n]\setminus A)$ with $\gamma n\le|A|\le n/2$ contains at least $\rho n^2$ edges of $H_1$. Conditional on such a realization, the probability that some such cut contains no edge of $H_{p_n^-}$ is at most
$2^n e^{-\rho p_n^-n^2}=o(1)$, since $np_n^-=C_n\to\infty$. As there are at most $1/\gamma$ components of order at least $\gamma n$, we obtain
\begin{equation}\label{eq:small-components-macro}
\EE\Big[\textstyle\sum_{\gamma n\le j\le n/2}\mathcal C_j^-\Big]=o(1).
\end{equation}
Combining~\eqref{eq:small-components-small-range}--\eqref{eq:small-components-macro} thus proves~\ref{lem:preconnectivity-window:itm2}.
\end{proof}
By \Cref{lem:preconnectivity-window}, the proof of \Cref{theorem:hitting-time} now follows the same structure as the classical Erd\H{o}s--R\'enyi argument: once the graph consists of a single non-trivial component and isolated vertices, it suffices to show that no edge joining two of these isolated vertices appears before the graph becomes connected.

\begin{proof}[Proof of Theorem~\ref{theorem:hitting-time}] Let $G\sim G(n,W)$. Take $p_n^-,p_n^+$ and $h_n$ from \Cref{lem:preconnectivity-window}, and let
$K_-:=Y^{(p_n^-)}$ and~$K_+:=Y^{(p_n^+)}$  denote the numbers of edges of $H_{p_n^-}$ and $H_{p_n^+}$, respectively. {As mentioned at the beginning of \Cref{sec:dynamic}}, we have that~$H_{p_n^-}=G_{K_-}$ and $H_{p_n^+}=G_{K_+}$.

Let $V_1$ denote the set of isolated vertices of $H_{p_n^-}$ and observe that, with probability tending to one, \Cref{lem:preconnectivity-window} implies that $1\le |V_1|\le2h_n$ and that $H_{p_n^-}$ has no component of order between $2$ and $n/2$. Hence, conditioned on this event, $H_{p_n^-}$ consists of one connected component together with $V_1$. The same lemma gives that $H_{p_n^+}$ is connected with probability tending to one.

We claim that, with probability tending to one, no edge born during the interval $(p_n^-,p_n^+]$ has both endpoints in $V_1$. 
For a fixed pair $u,v$, let $\mathcal F_{uv}$ be the
$\sigma$-algebra generated by all vertex types
$X_1,\ldots,X_n$ and all priority variables $U_{ab}$
with $ab\ne uv$. The absence of all edges incident to
$u$ or $v$, other than $uv$, in $H_{p_n^-}$ is
$\mathcal F_{uv}$-measurable. Thus, conditional on
$\mathcal F_{uv}$ and on $u,v\in V_1$, the only
restriction on $U_{uv}$ is
$U_{uv}>p_n^-W(X_u,X_v)$. Since $U_{uv}$ is independent
of $\mathcal F_{uv}$, we obtain
\[
\PP\big(uv\in E(H_{p_n^+})\setminus E(H_{p_n^-})
\mid \mathcal F_{uv},\ u,v\in V_1\big)
\le
\frac{(p_n^+-p_n^-)W(X_u,X_v)}{1-p_n^-W(X_u,X_v)}
\le2p_n^+
\]
for all sufficiently large $n$, since $W\le1$ and $p_n^-\le1/2$. Summing this estimate over the fixed pairs $u<v$ gives
\[
\PP\big(\text{there is a } uv\in E(H_{p_n^+})\setminus E(H_{p_n^-})\text{ such that } u,v\in V_1\big)
\le 2p_n^+\EE\Big[\binom{|V_1|}{2}\Big].
\]
Since $|V_1|=\mathcal C_1^-$ and the proof of \Cref{lem:preconnectivity-window} gives
$\EE[(\mathcal C_1^-)_2]=(1+o(1))h_n^2$, the right-hand side is~$O(p_n^+h_n^2)=o(1)$,
where the last equality follows from
\Cref{lem:preconnectivity-window}.

Since a.a.s.\ no edge appearing between $p_n^-$ and $p_n^+$ joins two vertices of $V_1$, and $H_{p_n^-}$ has a unique non-trivial component, each vertex of $V_1$ loses its isolation by acquiring an edge to that component. Consequently, throughout the interval from $K_-$ to $K_+$ the graph consists of one non-trivial connected component and a subset of $V_1$ of isolated vertices. Since at least one isolated vertex is present at step $K_-$ and the graph is connected at step $K_+$, the moment at which the last isolated vertex disappears is exactly the step at which the graph becomes connected and therefore $\tau_{\mathrm{conn}}(W)=\tau_1(W)$.
\end{proof}

\subsection{Location and sharpness of the connectivity transition}\label[section]{sec:location}
In this section, we give a simple analytic condition, expressed in terms of the Laplace transform $\mathcal L(C)=I_0(C)$, under which, for a connected graphon $W$ satisfying $g_W(t)/t\to0$ as $t\to0$, the hitting times $\tau_1(W)$ and~$\tau_{\mathrm{conn}}(W)$ lie a.a.s.\ in a common deterministic window whose endpoints have ratio tending to one. To do so, let us suppose that there are sequences $C_n\to\infty$ and $\varepsilon_n\to0$ such that
\begin{equation}\label{eq:sharp-laplace-window}
\lim_{n\to\infty}n\mathcal L\big((1-\varepsilon_n)C_n\big)=\infty
\quad\text{and}\quad
\lim_{n\to\infty}n\mathcal L\big((1+\varepsilon_n)C_n\big)=0.
\end{equation}

\begin{proposition}\label[proposition]{prop:sharp-window} Let $W$ be a connected graphon such that $\lim_{t\to 0}\frac{g_W(t)}t=0$ and~\eqref{eq:sharp-laplace-window} holds. Then, there exist deterministic sequences $k_{-}=k_{-}(n)$ and $k_{+}=k_{+}(n)$  such that $k_-<k_+$,
\[
\frac{k_-}{k_+}\xrightarrow[n\to\infty]{}1
\qquad
\text{and}
\qquad
\PP\big(k_-<\tau_1(W)=\tau_{\mathrm{conn}}(W)\le k_+\big)\xrightarrow[n\to\infty]{}1.
\]
\end{proposition}

\begin{proof}
Note that we may enlarge $\varepsilon_n$ while keeping $\varepsilon_n\to0$, as both limits in~\eqref{eq:sharp-laplace-window} are preserved by monotonicity of $\mathcal L$. Thus we may assume, say, that $\varepsilon_n\ge n^{-1/4}$.

We first note that $C_n=o(n)$. Indeed, with $C_n^-:=(1-\varepsilon_n)C_n$, \Cref{lem:Ir-small-tail} gives $C_n^- \mathcal L(C_n^-)\to0$, whereas $n \mathcal L(C_n^-)\to\infty$. Hence $C_n^-/n\to0$, and therefore $C_n/n\to0$. Let $a_n:=n\mathcal L\big((1-\varepsilon_n)C_n\big)$ and $h_n:=\min\{a_n^{1/2},(\log n)^{1/2}\}$.
Then $h_n\to\infty$ and $h_n\le\log n$. Since $\mathcal L$ is continuous and decreasing, for all sufficiently large $n$ there is a unique $\widehat C_n$ with $(1-\varepsilon_n)C_n\le \widehat C_n\le (1+\varepsilon_n)C_n$ such that
$n\mathcal L(\widehat C_n)=h_n$.
In particular, $\widehat C_n/C_n\to1$. Put
$\widehat p_n:=\widehat C_n/n$. Recall that $N^{(0)}_p$ is the number of isolated vertices in $H_p$. Applying \Cref{lem:thinned-isolated-factorial} with $j=1,2$ then gives
$\EE[N^{(0)}_{\widehat p_n}]=(1+o(1))h_n$ and
$\EE[(N^{(0)}_{\widehat p_n})_2]=(1+o(1))h_n^2$, so~$\operatorname{Var}(N^{(0)}_{\widehat p_n})=o(h_n^2)$. Since $h_n\to\infty$, Chebyshev's inequality yields
$\PP(\delta(H_{\widehat p_n})=0)\longrightarrow 1$.

To obtain an upper bound on the time at which the last isolated vertex disappears, set $p_n^+:=(1+2\varepsilon_n)C_n/n$. Since $\varepsilon_n\ge n^{-1/4}$, for all sufficiently large $n$ we have
$(n-1)p_n^+\ge(1+\varepsilon_n)C_n$. Thus, the expected number of isolated vertices satisfies
\[
\EE[N^{(0)}_{p_n^+}]
\le n\mathcal L((n-1)p_n^+)
\le n\mathcal L\big((1+\varepsilon_n)C_n\big)
=o(1).
\]
Hence, by Markov's inequality, $H_{p_n^+}$ has no isolated vertex a.a.s.

We now translate the two thinning times to deterministic edge counts. Recall that $Y^{(p)}=e(H_p)$. Conditional on $G=H_1$, we have $Y^{(p)}\sim\operatorname{Bin}(Y^{(1)},p)$, and hence
$\EE\big[(Y^{(p)}-pY^{(1)})^2\big]=p(1-p)M_n$.
Since $Y^{(1)}/M_n\to1$ in probability by \Cref{lem:Yconc}
and $pM_n\to\infty$ whenever $np\to\infty$, it follows
that $Y^{(p)}/(pM_n)\to1$ in probability. Choose any sequence~$(\delta_n)_{n\in\mathbb N}$ such that $\delta_n\to0$ sufficiently slowly so that
\begin{equation}\label{eq:sharp-edge-count-concentration}
\PP\Big(
Y^{(\widehat p_n)}>(1-\delta_n)\widehat p_nM_n,
\
Y^{(p_n^+)}<(1+\delta_n)p_n^+M_n
\Big)\xrightarrow[n\to\infty]{}1.
\end{equation}
Set
$k_-:=\lfloor(1-\delta_n)\widehat p_nM_n\rfloor$ and
$k_+:=\lceil(1+\delta_n)p_n^+M_n\rceil$.
Conditioned on the event given by \eqref{eq:sharp-edge-count-concentration}, $G_{k_-}\subseteq H_{\widehat p_n}$ and $H_{p_n^+}\subseteq G_{k_+}$. Thus a.a.s.\ $G_{k_-}$ has an isolated vertex and $G_{k_+}$ has no isolated vertices. It follows that
$k_-<\tau_1(W)\le k_+$ with probability tending to one, and, moreover,
\[
\frac{k_-}{k_+}
=\big(1+o(1)\big)
\frac{\widehat C_n}{(1+2\varepsilon_n)C_n}
\longrightarrow1.
\]
The conclusion thus follows from \Cref{theorem:hitting-time}.
\end{proof}
We now set $d_*:=\operatorname*{ess\,inf}_{x\in[0,1]}d_W(x)$. For each fixed $\ell\ge1$, the next result characterizes when there is a deterministic time $k=\Theta(n\log n)$ at which $G_k$ has minimum degree at least $\ell$ a.a.s.
\begin{proposition}\label[proposition]{prop:nlogn-essential-inf}
Let $W$ be a graphon with $\|W\|_1>0$ and let $\ell\ge1$. There exists a sequence $k=k(n)=\Theta(n\log n)$ such that $\PP(\tau_\ell(W)\le k)\to1$ if and only if $d_*>0$. 
\end{proposition}
\begin{proof}Suppose first that $d_*>0$ and choose $\delta\in (0,d_*)$. Since $M_n=(1+o(1))\|W\|_1n^2/2$, if $k=an\log n$ for some fixed $a>0$, then $C_n:=nk/M_n=(1+o(1))(2a/\|W\|_1)\log n$. Moreover, by definition, $d_W(x)\ge\delta$ almost everywhere, and therefore, for every fixed $r\ge0$, we have
$I_r(C_n/2)\le e^{-\delta C_n/2}$. Hence
\[
nC_n^r I_r(C_n/2)
\le
nC_n^r e^{-\delta C_n/2}
=
n^{1-a\delta/\|W\|_1+o(1)}(\log n)^r.
\]
Choosing $a>\|W\|_1/\delta$ makes the right-hand side tend to zero for every fixed $r$. Since $k\rightarrow\infty$ and $k = o(M_n)$, \Cref{cor:Z-parameters} applies and, together with \Cref{prop:Z-two-sided}, yield
\begin{equation*}\EE[N_k^{(r)}] \leq n^{1-a\delta/\|W\|_1+o(1)}(\log n)^r = o(1)
\end{equation*}
for every fixed $r\geq0$. Therefore, by Markov's inequality $\PP(N_k^{(\leq\ell)}\geq 1) =o(1)$, and thus~$\PP(\tau_\ell(W)\le k)\to1$.

For the other direction, we will prove that contrapositive. So let us assume that $d_*=0$, and let $k=k(n)=\Theta(n\log n)$. Fix $\xi>0$ and use the thinning coupling from the beginning of \Cref{sec:dynamic} with $p_+:=(1+\xi)k/M_n$. Then $p_+=O(\log n/n)$, and for some positive constant~$A$, we have $np_+\le A\log n$ for all sufficiently large $n$. Choose $\varepsilon>0$  small enough so that $A\varepsilon<1/2$, and put~$S:=\{x:d_W(x)\le\varepsilon\}$. Since $d_*=0$, we have that $s:=\lambda(S)=g_W(\varepsilon)$ is positive.
Let $H:= H_{p_+} \sim G(n, p_+ W)$ and let $V_n$ count isolated vertices of $H$ whose types lie in $S$. Observe that 
\[
\EE[V_n] = n \int_S (1 - p_+ d_W(x))^{n-1} \, dx.
\]
Defining $J_n:= \int_S e^{-np_+ d_W(x)} \, dx$, the condition $np_+^2 = o(1)$ implies that $\EE[V_n] = (1 + o(1)) n J_n.$ Since $d_W(x) \le \varepsilon$ for $x \in S$, it follows that $J_n \ge s e^{-np_+\varepsilon} \ge s n^{-A\varepsilon}$, and therefore $\EE[V_n] \to \infty$. For distinct vertices with types $x,y\in S$, the probability that both are isolated in $H$ equals
\[
(1-p_+W(x,y))
\big(1-p_+d_W(x)-p_+d_W(y)+p_+^2q_W(x,y)\big)^{n-2}.
\]
Since $0\le q_W(x,y)\le1$ and $np_+^2=o(1)$, this is at most
$(1+o(1))e^{-np_+(d_W(x)+d_W(y))}$ uniformly on $S^2$. It follows that
$\EE[(V_n)_2]\le(1+o(1))n^2J_n^2=(1+o(1))\EE[V_n]^2$. Since $\EE[V_n]\to\infty$, the second-moment inequality yields $\PP(V_n>0)\to1$.

Then, with probability tending to one, the thinning coupling gives $G_k\subseteq H$, and thus every isolated vertex in $H$ is also isolated in $G_k$. Hence $\PP(\delta(G_k)=0)\to1$, or equivalently $\PP(\tau_1(W)>k)\to1$.
\end{proof}

\begin{proposition}\label{prop:precise hitting time}Let $W$ be a connected graphon with $d_*>0$. Then,   $\frac{\tau_1(W)}{n\log n}\xrightarrow{\mathbb P}\frac{\|W\|_1}{2d_*}$ and $\frac{\tau_{\mathrm{conn}}(W)}{n\log n}\xrightarrow{\mathbb P}\frac{\|W\|_1}{2d_*}$. In particular, if $W\equiv p$ for $p\in (0,1]$, then $
\mathbb P\big(\tau_1(W)=\tau_{\mathrm{conn}}(W)\big)\to1$ and
$\frac{\tau_1(W)}{n\log n}\xrightarrow{\mathbb P}\frac12$.\end{proposition}
\begin{proof}Since $d_W(x)\ge d_*$ almost everywhere, we have $g_W(t)=0$ for $t<d_*$ and thus $\lim_{t\to 0}\frac{g_W(t)}{t}=0$. Therefore, by \Cref{theorem:hitting-time}, we have  a.a.s.\ that $\tau_1(W)=\tau_{\mathrm{conn}}(W)<+\infty$. Now, note that 
\[\mathcal L(C)=\int_0^1e^{-Cd_W(x)}dx\le e^{-Cd_*}\]
and, as for a fixed $\varepsilon>0$ we have $g_W(d_*+\varepsilon)>0$, we have $\mathcal L(C)\ge e^{-C(d_*+\varepsilon)}g_W(d_*+\varepsilon)$. Therefore, 
\[d_*\le \frac{-\log\mathcal L(C)}{C}\le d_*+\varepsilon-\frac{\log g_W(d_*+\varepsilon)}{C},\]
from which follows that $\lim_{C\to +\infty}\frac{-\log\mathcal L(C)}{C}=d_*$. For $n\in\mathbb N$ choose $C_n$ so that $n\mathcal L(C_n)=1$, and note that $C_n=(1+o(1))\frac{\log n}{d_*}$. Let $\varepsilon_n:=(\log n)^{-1/2}$ and note that
\[\mathcal L((1+\varepsilon_n)C_n)\le e^{-\varepsilon_n C_nd_*}\mathcal L(C_n)\quad \text{and}\quad \mathcal L((1-\varepsilon_n)C_n)\ge e^{\varepsilon_nC_nd_*}\mathcal L(C_n),\]
from which \eqref{eq:sharp-laplace-window} easily follows. From the proof of \Cref{prop:sharp-window}, it follows that a.a.s.\ $\tau_1(W)=\tau_{\mathrm{conn}}(W)=(1+o(1))C_nM_n/n$, where $M_n=\binom{n}{2}\|W\|_1$. Using that $C_n=(1+o(1))\frac{\log n}{d_*}$, we conclude the first part of the proof. For the second part, we just need to note that $d_*=p$ and $\|W\|_1=p$ when $W\equiv p$. 
\end{proof}
The proof of Proposition~\ref{prop:precise hitting time} shows that condition~\eqref{eq:sharp-laplace-window} holds whenever the graphon degree function is bounded away from zero. Indeed, if~$d_*=\operatorname*{ess\,inf}d_W>0$, choosing $C_n$ so that $n\mathcal L(C_n)=1$ and setting $\varepsilon_n:=(\log n)^{-1/2}$ gives the two required limits. We now show that the condition can also hold when $d_*=0$. We demonstrate this with examples of graphons having stretched-exponential lower tails, for which the connectivity windows are multiplicatively sharp at scales larger than $n\log n$. 

\begin{example}{Fix $\alpha>0$ and put
\[
h_\alpha(x):=\big(\log(e/x)\big)^{-1/\alpha},
\qquad
\rho_\alpha:=\int_0^1h_\alpha(x)\,dx,
\]
and consider the  graphon $W_\alpha(x,y):=h_\alpha(x)h_\alpha(y)$.
Then $W_\alpha$ is connected and
$\operatorname*{ess\,inf}d_{W_\alpha}=0$. Moreover,
\[
d_{W_\alpha}(x)=\rho_\alpha\big(\log(e/x)\big)^{-1/\alpha}.
\]
Consequently, for $0<t\le\rho_\alpha$,
\[
g_{W_\alpha}(t)
=
\exp\Big(1-\Big(\frac{\rho_\alpha}{t}\Big)^\alpha\Big).
\]
In particular,
$\log g_{W_\alpha}(t)=-(1+o(1))\rho_\alpha^\alpha t^{-\alpha}$ as $t\to0$. The corresponding Laplace transform (using the change of variable $y=\log(e/x)$) is
\[
\mathcal L_\alpha(C):=\int_0^1 e^{-C d_{W_\alpha}(x)}\,dx
=e\int_1^\infty
\exp\Big(-y-\rho_\alpha C y^{-1/\alpha}\Big)\,dy.
\]
Put
$\beta:=\alpha/(\alpha+1)$ and
$y_C:=(\rho_\alpha C/\alpha)^\beta$. The function $y+\rho_\alpha C y^{-1/\alpha}$ attains its minimum at $y_C$, where its second derivative equals $(\alpha+1)/(\alpha y_C)$. A second-order Laplace expansion around this minimum gives
\[
\mathcal L_\alpha(C)
=(1+o(1))e\sqrt{\frac{2\pi\alpha}{\alpha+1}}\,
y_C^{1/2}
\exp\big(-(\alpha+1)y_C\big).
\]
Equivalently, if
$\kappa_\alpha:=(\alpha+1)(\rho_\alpha/\alpha)^\beta$, then
$\log \mathcal L_\alpha(C)=-\kappa_\alpha C^\beta+O(\log C)$.
Therefore, with
\[
C_n:=\Big(\frac{\log n}{\kappa_\alpha}\Big)^{1/\beta}
=\Big(\frac{\log n}{\kappa_\alpha}\Big)^{(\alpha+1)/\alpha},
\]
and $\varepsilon_n=(\log n)^{-1/2}$, say,
condition~\eqref{eq:sharp-laplace-window} holds. Since $\|W_\alpha\|_1=\rho_\alpha^2$, 
the relation $k=M_nC_n/n$ gives $M_nC_n/n=\Theta(n(\log n)^{1+1/\alpha})$.
Thus, for every $\alpha>0$, the connectivity transition is localized in a $(1+o(1))$-multiplicative window at a scale larger than $n\log n$, despite $\operatorname*{ess\,inf}d_{W_\alpha}=0$.

The same conclusion depends only on the logarithmic lower-tail asymptotic, not on the rank-one representation. If for some $a,\alpha>0$, $\log g_W(t)=-(1+o(1))a t^{-\alpha}$ as $t\to 0$,
then integration by parts gives~$\mathcal L(C)=e^{-C}+C\int_0^1e^{-Ct}g_W(t)\,dt$. After the scaling $t=C^{-1/(\alpha+1)}s$, the leading exponential behavior of the integral is determined by minimizing $s+as^{-\alpha}$ over $s> 0$. This gives
\[
\log \mathcal L(C)=-(1+o(1))\kappa C^{\alpha/(\alpha+1)},
\qquad
\kappa:=(\alpha+1)a^{1/(\alpha+1)}
\alpha^{-\alpha/(\alpha+1)}.
\]
Choosing
$C_n=(\log n/\kappa)^{(\alpha+1)/\alpha}$ and then choosing
$\varepsilon_n\to0$ sufficiently slowly that
$\varepsilon_n\log n\to\infty$ and the $o(1)$ error in the last logarithmic asymptotic is $o(\varepsilon_n)$ gives~\eqref{eq:sharp-laplace-window}. Thus stretched-exponential lower tails provide a broad class of graphons with zero essential infimum but a multiplicatively sharp connectivity window.}
\end{example}

\section{Discussion and open problems}\label[section]{sec:discussion}
We conclude with three open problems concerning hitting times for stronger graph properties in the \(W\)-random graph process.

\medskip
\noindent\textbf{Hitting times for $K$-connectivity. } As mentioned in the introduction,  Hladk\'y and Viswanathan~\cite{HV} asked whether the hitting time for connectivity and minimum degree at least $K$ coincide in their $W$-edge-incremental process. Given our \Cref{theorem:hitting-time}, it is natural to ask the same question for the $W$-random graph process. To state this question precisely, let $W$ be a graphon and let $(G_t)_{t=0}^{e(G(n,W))}$  be the $W$-random graph process, and, for~$K\ge 1$, let  $\tau_{K-\mathrm{conn}}(W):=\min\{t:G_t\text{ is $K$-connected}\}$ and $\tau_K(W):=\min\{t:\delta(G_t)\ge K\}$. 
\begin{conjecture}Let $W$ be a connected graphon such that $\lim_{t\to 0}\frac{g_W(t)}{t}=0$. Then, for every $K\ge 2$, 
\[\lim_{n\to\infty}\mathbb P\big(\tau_{K-\mathrm{conn}}(W)=\tau_K(W)\big)=1.\] 
\end{conjecture}

\noindent\textbf{Hamiltonicity hitting time. } The study of Hamiltonicity of random graphs has been one of the central problems in the area. Indeed, for the Erd\H{o}s--R\'enyi random graph process,
Bollob\'as~\cite{Bollobas1984} and, independently, Ajtai, Koml\'os and
Szemer\'edi~\cite{AjtaiKomlosSzemeredi1985} showed that a.a.s.~the process becomes Hamiltonian
exactly when its minimum degree reaches two. For graphons, we may ask the same questions. Let $W$ be a graphon and let~$(G_t)_{t=0}^{e(G(n,W))}$ be the $W$-random graph process. Define the hitting time for Hamiltonicity as
\[\tau_{\mathrm{Ham}}(W):=\min\{t:G_t\text{ is Hamiltonian}\},\]
where $\tau_{\mathrm{Ham}}(W)=+\infty$ whenever no $G_t$ is Hamiltonian.  Recently, Garbe, Hladk\'y and Piga~\cite{GarbeHladkyPiga} characterized the graphons $W$ for which $G(n,W)$ is a.a.s.~Hamiltonian. We say that a graphon $W$ is Hamiltonian if 
\begin{enumerate}[label=\upshape(\roman{enumi})]
    \item\label{ham:1} $W$ is connected,
    \item\label{ham:2}  $\lim_{t\to 0}\frac{g_W(t)}{t}=0$, and
    \item\label{ham:3} there is no $a\in (0,1/2]$ and measurable disjoint sets $A,B\subset [0,1]$ with $\lambda(A)=a$ and $\lambda(B)=1-2a$ such that~$W=0$ almost everywhere on $A\times (A\cup B)$.
\end{enumerate}
In~\cite{GarbeHladkyPiga} it is shown that for a Hamiltonian graphon $W$,  the graph $G(n,W)$ a.a.s.~contains a Hamilton cycle  and, moreover, conditions~\ref{ham:1}--\ref{ham:3} are necessary for Hamiltonicity.  It seems plausible that a graphon-analogue of the hitting time of Hamiltonicity should also be true. 
\begin{conjecture}Let $W$ be a graphon satisfying~\ref{ham:1}--\ref{ham:3}. Then, 
\[\lim_{n\to\infty}\mathbb P\big(\tau_{\mathrm{Ham}}(W)=\tau_2(W)\big)=1.\] 
\end{conjecture}


\noindent\textbf{Characterizing sharp connectivity transitions.}
A third direction concerns the sharp localization of the common connectivity/minimum-degree hitting time. \Cref{prop:sharp-window} provides a sufficient condition in terms of the Laplace transform $\mathcal L$, requiring the existence of sequences $C_n\to\infty$ and $\varepsilon_n\to 0$ satisfying \eqref{eq:sharp-laplace-window}. It would be interesting to identify more explicit sufficient conditions, and ideally a necessary and sufficient condition, for sharp localization in terms of the asymptotic decay of $\mathcal L(C)$ under multiplicative perturbations of $C$.

\section*{Acknowledgements}
The authors used ChatGPT (OpenAI) as an interactive assistant during the development and preparation of this manuscript, including for exploring proof strategies, checking mathematical arguments and calculations, and improving the exposition.

\appendix

\section{Appendix}\label[appendix]{app:static-simplifications}

For completeness, we collect here the proofs of \Cref{lem:exceptional} and \Cref{thm:static-connectivity}. Both results are due to Hladk\'y and Viswanathan~\cite{HV}, but we decided to include both arguments as we believe they give genuine simplifications of the corresponding proofs. The first replaces their auxiliary ``sprout'' technique by a direct double counting argument, while the second relies on an estimate involving the largest graphon degree in a candidate component to replace the dyadic calculations from~\cite{HV}. 

\subsection{The exceptional-set lemma}

\begin{proof}[Proof of \Cref{lem:exceptional}]
Let $S_0:=S_W(2c/n)$ and, recursively for $m\ge1$, define
\[
S_{m}:=\Big\{x\notin\bigcup_{i<m}S_i:
 d_W\Big(x,[0,1]\setminus\bigcup_{i<m}S_i\Big)<\frac14d_W(x)\Big\},
\qquad S:=\bigcup_{m\ge0}S_{m}.
\]
If $x\notin S$, then for every $r\ge1$ one has
$d_W(x,[0,1]\setminus\bigcup_{i<r}S_i)\ge d_W(x)/4$, otherwise $x$ would belong to $S_r$. The complements $[0,1]\setminus\bigcup_{i<r}S_i$ decrease to $[0,1]\setminus S$, so monotone convergence gives
$d_W(x,[0,1]\setminus S)=\lim_{r\to\infty}d_W(x,[0,1]\setminus\bigcup_{i<r}S_i)\ge d_W(x)/4$. Thus the degree-retention property follows from the construction. It remains to bound $\lambda(S)$.

For $0\le i<j$, set $w_{ij}:=\int_{S_i\times S_j}W$, and let
$T:=\sum_{i<j}w_{ij}$ and $R:=\sum_{j\ge1}w_{0j}$. If $m\ge1$, every $x\in S_{m}$ sends more than three quarters of its degree to the earlier layers, and hence
$\sum_{i<m}w_{im}\ge3\sum_{j>m}w_{mj}$. Summing over $m\ge1$ gives $T\ge 3(T-R)$, so $T\le3R/2$. Since $d_W(x)\le2c/n$ on $S_0$,
$R\le(2c/n)g_W(2c/n)$ and therefore $T\le(3c/n)g_W(2c/n)$.

If $x\in S_{m}$ with $m\ge1$, then $x\notin S_0$ and hence $d_W(x)>2c/n$. At least three quarters of this degree goes to earlier layers. 
Thus, for every $m\ge1$,
\[
\sum_{i<m}w_{im}
=\int_{S_m}\int_{\bigcup_{i<m}S_i}W(x,y)\,dy\,dx
\ge\frac{3c}{2n}\lambda(S_m).
\]
Summing over $m\ge1$ and using the fact that the sets $S_m$
are pairwise disjoint, we obtain
$T
\ge\frac{3c}{2n}
\lambda(\bigcup_{m\ge1}S_m)$.
Therefore, $\lambda(\bigcup_{m\ge1}S_{m})\le\frac{2n}{3c}T\le2g_W(2c/n)$. Adding $\lambda(S_0)=g_W(2c/n)$ proves $\lambda(S)\le3g_W(2c/n)$.
\end{proof}

\subsection{$K$-connectivity}

\begin{proof}[Proof of \Cref{thm:static-connectivity}]
We first prove that $G(n,W)$ is a.a.s.\ connected.

Fix $\varepsilon>0$. We shall choose $\gamma\in(0,1/2]$ so that, for all sufficiently large $n$,
the probability that $G(n,W)$ has a component of size at most $\gamma n$ is less than
$\varepsilon/2$. The possibility of a component of size between $\gamma n$ and~$n/2$ is then
ruled out by connectedness of $W$.

\medskip
\noindent
\textbf{Small components.}
For $A\subseteq[n]$, let
$D_A:=\{e_G(A,[n]\setminus A)=0\}$.
Fix $A$ with $|A|=s\le n/2$, and condition on the types
$(X_i)_{i\in A}=(x_i)_{i\in A}$. Set
$M_A:=\max_{i\in A}d_W(x_i)$.
If $i^\ast\in A$ attains this maximum, then
\begin{align*}
\PP(D_A\mid (X_i)_{i\in A}=(x_i)_{i\in A})
&=
\Big(
\int_0^1 \prod_{i\in A}(1-W(x_i,y))\,dy
\Big)^{n-s}\\
&\le
\Big(
\int_0^1(1-W(x_{i^\ast},y))\,dy
\Big)^{n-s}\\
&=(1-M_A)^{n-s}\\
&\le e^{-(n-s)M_A}.
\end{align*}
Unconditioning, and writing $\rho:=n-s$, we therefore obtain
$\PP(D_A)\le \EE[e^{-\rho M_A}]$.
Since the \(X_i\)'s are independent,
$\PP(M_A\le t)=g_W(t)^{s}$.
Moreover,
\[
e^{-\rho M_A}
=
\int_{M_A}^{\infty}\rho e^{-\rho t}\,dt,
\]
and so the Fubini--Tonelli theorem gives
\begin{equation}\label{eq:static-cut-laplace}
\PP(D_A)
\le
\rho\int_0^\infty e^{-\rho t}g_W(t)^{s}\,dt,
\end{equation}
where we set $g_W(t)=1$ for $t\ge1$.

Choose $\varepsilon_0>0$ so small that
$2\varepsilon_0/(1-2\varepsilon_0)<\varepsilon/4$.
By the assumption $\lim_{t\to 0} g_W(t)/t=0$, there exists~$\delta>0$ such that
$g_W(t)\le \varepsilon_0 t$
for every $t\in[0,\delta]$.
Splitting the integral in~\eqref{eq:static-cut-laplace} at $\delta$,
\begin{align*}
\PP(D_A)
\le
\rho\int_0^\delta e^{-\rho t}(\varepsilon_0 t)^{s}\,dt
+
\rho\int_\delta^\infty e^{-\rho t}\,dt
\le
\varepsilon_0^{s}\rho\int_0^\infty t^{s} e^{-\rho t}\,dt
+
e^{-\rho\delta}
&=
\varepsilon_0^{s}\frac{s!}{(n-s)^{s}}
+
e^{-(n-s)\delta}.
\end{align*}

We now choose \(\gamma\in(0,1/2]\) sufficiently small that
$H(\gamma)<(1-\gamma)\delta$,
where
$H(\gamma):=-\gamma\log\gamma-(1-\gamma)\log(1-\gamma)$
is the binary entropy function. A union bound and the preceding estimate yield
\begin{align*}
\PP\bigl(\exists A\subseteq[n]:
1\le |A|\le\gamma n,\ D_A\bigr)
&\le
\sum_{s=1}^{\lfloor\gamma n\rfloor}
\binom n{s}
\Big(
\varepsilon_0^{s}\frac{s!}{(n-s)^{s}}
+
e^{-(n-s)\delta}
\Big).
\end{align*}
The first contribution is at most
\[
\sum_{s=1}^{\lfloor\gamma n\rfloor}
\binom n{s}\varepsilon_0^{s}\frac{s!}{(n-s)^{s}}
\le
\sum_{s=1}^{\lfloor\gamma n\rfloor}
\Big(\frac{\varepsilon_0 n}{n-s}\Big)^{s}
\le
\sum_{s\ge1}(2\varepsilon_0)^{s}
=
\frac{2\varepsilon_0}{1-2\varepsilon_0}
<
\frac{\varepsilon}{4}.
\]
For the second contribution, the standard entropy bound gives
\[
\sum_{s=1}^{\lfloor\gamma n\rfloor}
\binom n{s} e^{-(n-s)\delta}
\le
\exp\Big(
-\bigl((1-\gamma)\delta-H(\gamma)-o(1)\bigr)n
\Big)
=o(1).
\]
Hence, for all sufficiently large \(n\),
\[
\PP\bigl(\exists A\subseteq[n]:
1\le |A|\le\gamma n,\ D_A\bigr)
<
\frac{\varepsilon}{2}.
\]

\medskip
\noindent
\textbf{Macroscopic components.}
For components of linear order, we use
\cite[Proposition~11]{HV}, which follows from Janson's
characterization of connected graph limits~\cite[Theorem~1.10(i)]{JansonConnected}.
Since $W$ is connected and $G(n,W)$ converges to $W$ in cut distance
a.a.s.~\cite[Lemma~10.16]{LovaszBook}, this proposition implies that, for every fixed $\gamma>0$,
there exists $\rho=\rho(W,\gamma)>0$ such that a.a.s.
$G(n,W)$ has at least $\rho n^2$ edges across every cut
$(A,[n]\setminus A)$ satisfying $\gamma n\le |A|\le n/2$.
Consequently,
\[
\PP\bigl(\exists A\subseteq[n]:
\gamma n\le |A|\le n/2,\ D_A\bigr)=o(1).
\]

\paragraph{Conclusion.}If $G(n,W)$ is disconnected, then one of its components has size at most
$n/2$. Combining the small- and macroscopic-component estimates therefore
shows that $G(n,W)$ is a.a.s.\ connected.

Hladk\'y and Viswanathan~\cite[Lemma~12]{HV} prove that, for any graphon
$W$, asymptotically almost sure connectivity of $G(n,W)$ implies asymptotically
almost sure $K$-connectivity for every fixed $K\in\N$. Applying their lemma
to the connectivity statement just proved completes the proof.
\end{proof}

\end{document}